\documentclass{article}

\usepackage{PRIMEarxiv}
\usepackage{multirow}
\usepackage[utf8]{inputenc} % allow utf-8 input
\usepackage[T1]{fontenc}    % use 8-bit T1 fonts
\usepackage{amsthm} % for theorem-like environments such as remark
\newtheorem{remark}{Remark}
\usepackage{hyperref}       % hyperlinks
\usepackage{url}            % simple URL typesetting
\usepackage{booktabs}       % professional-quality tables
\usepackage{amsfonts}       % blackboard math symbols
\usepackage{nicefrac}       % compact symbols for 1/2, etc.
\usepackage{microtype}      % microtypography
\usepackage{amsmath}        % for \boldsymbol
\usepackage{amssymb}
\usepackage{lipsum}
\usepackage{fancyhdr}       % header
\usepackage{graphicx}       % graphics
\graphicspath{{media/}}     % organize your images and other figures under media/ folder
\usepackage{subcaption}
\usepackage{stmaryrd}
\newcommand{\jump}[1]{\llbracket #1 \rrbracket}
\title{A Localized Hybrid Finite Element--Random Feature Method for Elliptic Interface Problems on Unfitted Meshes
}

\author{
  Siyuan Lang \\
  Ministry of Education Key Laboratory for NSLSCS, \\
  School of Mathematical Sciences, Nanjing Normal University \\
  \And
  Zhiyue Zhang\thanks{Corresponding author: \texttt{zhangzhiyue@njnu.edu.cn}} \\
  Ministry of Education Key Laboratory for NSLSCS, \\
  School of Mathematical Sciences, Nanjing Normal University \\
}

\begin{document}
\maketitle

\begin{abstract}
We develop a localized hybrid finite element--random feature method
(FEM--RFM) for elliptic interface problems on unfitted background meshes.
The central idea is to use the random feature method (RFM) only in the narrow band
formed by interface-cut cells of a uniform Cartesian background mesh,
while retaining standard $Q_p$ finite elements on the remaining uncut
cells. Within the band, local RFM representations are used on the
two physical sides of the interface. The prescribed jumps in the solution and normal flux are incorporated through residual equations involving the two RFM representations on the physical interface $\Gamma$, while the
finite element method (FEM) and RFM components are coupled across a separate grid-aligned
artificial interface $\Gamma_c$ through trace and variational normal-flux
residuals. A Schur reduction eliminates the interior finite element unknowns
through sparse solves while retaining their Galerkin equations as
equality constraints. The remaining band and interface residuals are
minimized in a reduced least-squares problem involving only the
coupling nodal values and random-feature coefficients. The method is assessed numerically on a range of two- and three-dimensional interface problems involving high coefficient contrasts, nonsmooth and multiple interfaces, and irregular outer boundaries, together with an application-oriented thermal-homogenization problem for a periodic bicontinuous composite.
\end{abstract}

% keywords can be removed

\section{Introduction}

Interface problems arise naturally in the mathematical modeling of
heterogeneous and multiphase media, including composite materials,
heat and mass transfer across material interfaces, and other processes
involving discontinuous physical properties
\cite{greengard1994numerical,chen1998finite}.
A representative elliptic interface problem involves a diffusion
coefficient that is discontinuous across an interface $\Gamma$,
together with prescribed jumps of the solution and its normal flux.
Although the solution may be sufficiently regular within each
subdomain, it is generally only piecewise regular across $\Gamma$. From the numerical viewpoint, the interface treatment must account for
the interface geometry, discontinuous coefficients, and prescribed jump conditions \cite{babuvska1970finite,leveque1994immersed,chen1998finite}.

Among classical mesh-based discretizations, a common distinction is
between interface-fitted and unfitted approaches. In interface-fitted
methods, the mesh is aligned with $\Gamma$ so that the interface is
represented by element boundaries. Standard conforming finite element
methods (FEMs) on fitted meshes have been studied extensively
\cite{barrett1987fitted,bramble1996finite,chen1998finite}.
High-order fitted finite element formulations, including isoparametric
approaches, have also been developed
\cite{li2010optimal,varsakelis2014numerical}.
Discontinuous Galerkin (DG), hybridizable discontinuous Galerkin (HDG), and weak
Galerkin (WG) formulations have also been developed on fitted meshes
\cite{cai2011discontinuous,huynh2013high,
huang2020high,mu2013weak}.
Interface-fitted virtual element methods, together with dedicated mesh
generation in two and three dimensions, have also been considered
\cite{chen2017interface}.
This geometric conformity provides a direct representation of the
interface, but mesh generation can become more involved for complicated
or evolving interfaces and may require remeshing.
Unfitted methods retain a background mesh that is independent
of $\Gamma$ and incorporate the interface treatment into the discrete
formulation. Representative examples include the immersed interface method
\cite{leveque1994immersed,li1998fast,adams2002immersed};
Cartesian-grid methods
\cite{liu2000boundary,li2003new,zhou2006high,bochkov2020solving};
immersed finite element (IFE) methods
\cite{li2004immersed,gong2008immersed,lin2015partially,
guo2020immersed};
immersed finite-volume element (IFVE) methods
\cite{ewing1999immersed,wang2021new};
unfitted Nitsche and cut finite element method (CutFEM) formulations
\cite{hansbo2002unfitted,burman2015cutfem};
and generalized, extended, and high-order unfitted finite element
methods
\cite{babuvska2017strongly,xiao2020high,
lehrenfeld2018analysis,chen2021adaptive,burman2026unfitted}.
In IFE methods, interface information is incorporated
into modified local approximation spaces on cut elements
\cite{li2004immersed,gong2008immersed,
guo2019higher,adjerid2024high}. Unfitted Nitsche and CutFEM
formulations treat the interface conditions variationally on
cut cells, often together with stabilization or related techniques for
small-cut configurations
\cite{hansbo2002unfitted,burman2015cutfem,
lehrenfeld2018analysis,burman2021unfitted}. For high-order schemes,
the interface geometry and interface-dependent discrete terms must be
treated consistently with the bulk approximation order
\cite{guo2019higher,lehrenfeld2018analysis,
chen2021adaptive,adjerid2024high}.
A mature literature provides rigorous stability and convergence theory
for these approaches. In these methods, the interface-specific
discretization is concentrated primarily in the neighborhood of the
physical interface.

In recent years, learning-based methods have also been developed for the
numerical solution of partial differential equations. Representative
examples include the deep Galerkin method (DGM) and physics-informed neural
networks (PINNs), in which network parameters are determined by minimizing
objectives constructed from the governing-equation and boundary residuals
\cite{sirignano2018dgm,raissi2019physics}.
For interface problems, a globally defined representation may need to accommodate piecewise regularity and possible jumps across the physical interface.
This has motivated interface-aware constructions based on separate or
piecewise neural representations on neighboring subdomains
\cite{wang2020mesh,he2022mesh,wu2022inn},
architectures designed to encode discontinuities or derivative
irregularities
\cite{hu2022discontinuity,tseng2023cusp,
fan2025novel,zheng2025ig,hu2026categorical},
and operator-learning formulations for parametric interface problems
\cite{wu2024solving,bi2025xi}.

Within this broader development, randomized approximation methods are
particularly relevant to the present work. Their hidden or feature
parameters are randomly generated and subsequently fixed, while only the
output coefficients remain unknown. This fixed-feature principle underlies the random feature method (RFM)~\cite{chen2022bridging}, extreme learning machines
(ELMs)~\cite{huang2006extreme}, and randomized neural network formulations,
and avoids repeated nonlinear optimization of all network parameters. Randomized trial spaces have been combined
with collocation, domain decomposition, Petrov--Galerkin, and DG
formulations for general PDEs
\cite{calabro2021extreme,dong2021local,
shang2023randomized,sun2024local}.
For interface problems, recent developments include RFM, piecewise ELM,
local randomized neural network (LRaNN), discontinuity-capturing RFM, and
discontinuous or domain-separated ELM formulations
\cite{chi2024random,liang2025piecewise,li2025local,
lin2025discontinuous,song2026discontinuity,zhu2026extreme}.
For broad classes of linear interface equations, these fixed-feature representations lead to problems that are linear in the output coefficients.

Neural and randomized approximations have also been combined with
classical numerical discretizations
\cite{mitusch2021hybrid,fan2026hybrid,
wu2026dfenn,YANG2026}. One direction employs neural or
randomized trial spaces within conventional variational frameworks,
including Petrov--Galerkin, DG, and hybridized
Petrov--Galerkin formulations
\cite{shang2023randomized,sun2024local,dang2024local}.
Another enriches a finite element trial space with network-generated
functions, so that the FEM and neural components contribute
within the same physical region
\cite{chen2026pinn,guo2026neural}.
For interface problems, the numerical treatment of the interface
neighborhood has a specific role: representing the piecewise solution
and incorporating the prescribed jumps on an unfitted mesh. Cells not
intersected by $\Gamma$ lie entirely within individual material
subdomains and admit the standard FEM treatment. This motivates an interface-induced decomposition of approximation methods, retaining the standard FEM discretization and its sparse interior solves on uncut cells while using a localized random-feature representation near the interface.

Motivated by this viewpoint, we develop a localized hybrid finite
element--random feature method (FEM--RFM) on unfitted background meshes.
We first construct a uniform Cartesian background mesh independently of the position and shape of the physical interface. The cells intersected by the interface are then collected into a narrow band for the random-feature approximation, while standard \(Q_p\) finite elements are retained on the remaining uncut cells. Within the interface band, the solution on each side of $\Gamma$ is
represented by an independent RFM approximation, implemented using a
local extreme learning machine (local ELM) with fixed hidden parameters. The prescribed jumps of solution
and normal flux are imposed directly on $\Gamma$, whereas the FEM
and RFM approximations are coupled across a separate
grid-aligned artificial interface $\Gamma_c$ through a nodal trace
residual and a variational normal-flux residual.
The interior FEM unknowns are eliminated by a Schur reduction using
sparse linear solves, with their Galerkin equations retained as
equality constraints. Only the coupling nodal values on $\Gamma_c$
and the random-feature coefficients remain in the reduced
least-squares problem. The FEM interior coefficients are recovered
from the coupling values after the reduced problem has been solved.
The formulation therefore separates the two approximations both
spatially and algebraically: the band equations and interface conditions
are treated through residual minimization, whereas the interior FEM
equations are retained as equality constraints.
The resulting computational cost is governed by the interaction
between the standard finite element discretization and the reduced
coupled system, rather than by the nominal number of approximation
coefficients alone. High-order bulk approximation uses standard
$Q_p$ spaces on the uncut cells. This separates the construction of
the bulk approximation from the treatment of the physical interface,
although the attainable accuracy still depends on adequate band
approximation, residual sampling, and coupling quadrature. Numerical experiments examine high coefficient contrasts, nonsmooth and
multiple interfaces, irregular outer boundaries, three-dimensional
geometries, and high-order $Q_p$ discretizations. A periodic bicontinuous
composite is further considered as an application-oriented
thermal-homogenization example. Finally, quasilinear elliptic and
parabolic moving-interface extensions are investigated numerically.

The remainder of this paper is organized as follows.
Section~\ref{sec:problem} introduces the elliptic interface problem and
the notation used throughout the paper.
Section~\ref{sec:methodology} presents the hybrid FEM--RFM formulation,
including the interface band construction, the FEM and
random-feature approximation spaces, the FEM weak formulation, the Schur
reduction, and the coupling and physical interface residuals.
Section~\ref{sec:numerical_examples} reports numerical experiments for
two- and three-dimensional interface problems and discusses the observed
accuracy, convergence behavior, and several extensions of the proposed
framework. Finally, the main findings and directions for further
development are summarized in the concluding section.

\section{Problem formulation}
\label{sec:problem}

Let $d\geq2$, and let
$\Omega\subset\mathbb{R}^{d}$ be a bounded Lipschitz domain.
Let $\Gamma\subset\Omega$ be a Lipschitz interface that separates
$\Omega$ into two disjoint open subdomains $\Omega^-$ and $\Omega^+$,
so that
\begin{equation}
\label{eq:domain-decomposition}
    \Omega
    =
    \Omega^-\cup\Gamma\cup\Omega^+,
    \qquad
    \Omega^-\cap\Omega^+=\emptyset,
    \qquad
    \Gamma
    =
    \partial\Omega^-\cap\partial\Omega^+ .
\end{equation}
The unit normal $\boldsymbol n$ on $\Gamma$ is oriented from
$\Omega^-$ to $\Omega^+$.
The diffusion coefficient is defined piecewise by
\begin{equation}
\label{eq:piecewise-beta}
    \beta(\boldsymbol{x})
    =
    \begin{cases}
        \beta^-(\boldsymbol{x}),
        &\boldsymbol{x}\in\Omega^-,
        \\[1mm]
        \beta^+(\boldsymbol{x}),
        &\boldsymbol{x}\in\Omega^+ ,
    \end{cases}
\end{equation}
and is assumed to satisfy
\begin{equation}
\label{eq:beta-assumption}
    \beta^\pm\in W^{1,\infty}(\Omega^\pm),
    \qquad
    0<\beta_{\min}
    \leq
    \beta^\pm(\boldsymbol{x})
    \leq
    \beta_{\max}
    <\infty .
\end{equation}
Piecewise constant diffusion coefficients are included in this setting.
For a quantity $q$ defined piecewise on the two subdomains, let
$q^\pm$ denote its restrictions to $\Omega^\pm$ and the corresponding
one-sided traces on $\Gamma$. We define the jump across $\Gamma$ by
\begin{equation}
\label{eq:jump-definition}
    \jump{q}
    =
    q^+|_{\Gamma}
    -
    q^-|_{\Gamma}.
\end{equation}

For $s\in\{-,+\}$, we write
\[
    \partial_{\boldsymbol n}u^s
    =
    \nabla u^s\cdot\boldsymbol n .
\]
Accordingly,
\begin{equation}
\label{eq:flux-jump-definition}
    \jump{\beta\partial_{\boldsymbol n}u}
    =
    \beta^+
    \nabla u^+\cdot\boldsymbol n
    -
    \beta^-
    \nabla u^-\cdot\boldsymbol n .
\end{equation}

Given a source term $f$, Dirichlet boundary data $g$, and prescribed
interface jumps $w$ and $v$, we consider the elliptic interface problem
\begin{equation}
\label{eq:linear-interface-problem}
\begin{cases}
-\nabla\cdot
\bigl(
\beta(\boldsymbol{x})\nabla u
\bigr)
+\kappa u
=
f,
&
\boldsymbol{x}\in\Omega^-\cup\Omega^+,
\\[1mm]
\jump{u}
=
w,
&
\boldsymbol{x}\in\Gamma,
\\[1mm]
\jump{\beta\partial_{\boldsymbol n}u}
=
v,
&
\boldsymbol{x}\in\Gamma,
\\[1mm]
u
=
g,
&
\boldsymbol{x}\in\partial\Omega ,
\end{cases}
\end{equation}
where $\kappa\geq0$ is a constant.
For nonhomogeneous interface jumps, the solution is naturally
interpreted piecewise with respect to the material subdomains.
In particular, a nonzero solution jump permits a discontinuity of
$u$ across $\Gamma$, so global smoothness across the interface cannot
in general be assumed. We therefore work with the subdomain
restrictions $u^\pm$ and interpret regularity piecewise in
$\Omega^-$ and $\Omega^+$; see, for example,
\cite{gong2008immersed}.
\section{Methodology}
\label{sec:methodology}

We discretize problem~\eqref{eq:linear-interface-problem} using a
hybrid FEM--RFM formulation. Within the interface band, the two
one-sided RFM approximations are realized by local extreme learning
machines (local ELMs) with fixed hidden parameters. On the prescribed
background mesh, the cells intersected by the physical interface define
the random-feature band, and the remaining uncut cells form the FEM
region.
With the features fixed, the band approximations depend linearly
on the output coefficients. For the linear interface model, eliminating the
interior FEM unknowns yields a linear least-squares problem for
the remaining unknown coefficients.

For clarity, the formulation is first presented with a single local
ELM representation on each physical side of one interface band, with
independent representations used on the two physical sides.
The FEM and random-feature representations are coupled across the
grid-aligned artificial interface $\Gamma_c$ through a nodal trace
residual and a variational normal-flux residual, whereas the
prescribed solution and flux jumps are imposed separately on the
physical interface $\Gamma$.

\begin{figure}[t]
    \centering
    \begin{subfigure}[t]{0.49\linewidth}
        \centering
        \includegraphics[width=\linewidth]{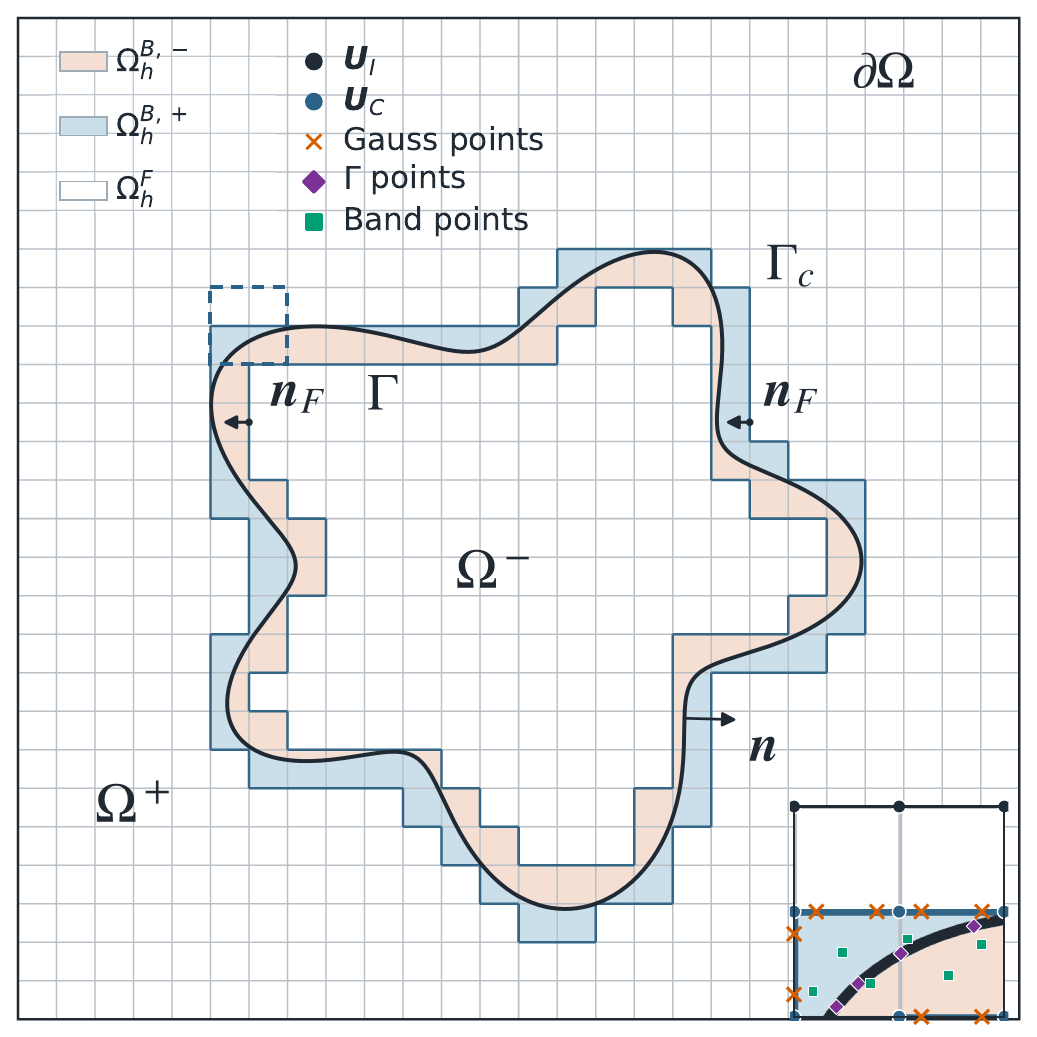}
        \caption{Hybrid domain decomposition.}
        \label{fig:hybrid_domain_decomposition}
    \end{subfigure}
    \hfill
    \begin{subfigure}[t]{0.49\linewidth}
        \centering
        \includegraphics[width=\linewidth]
        {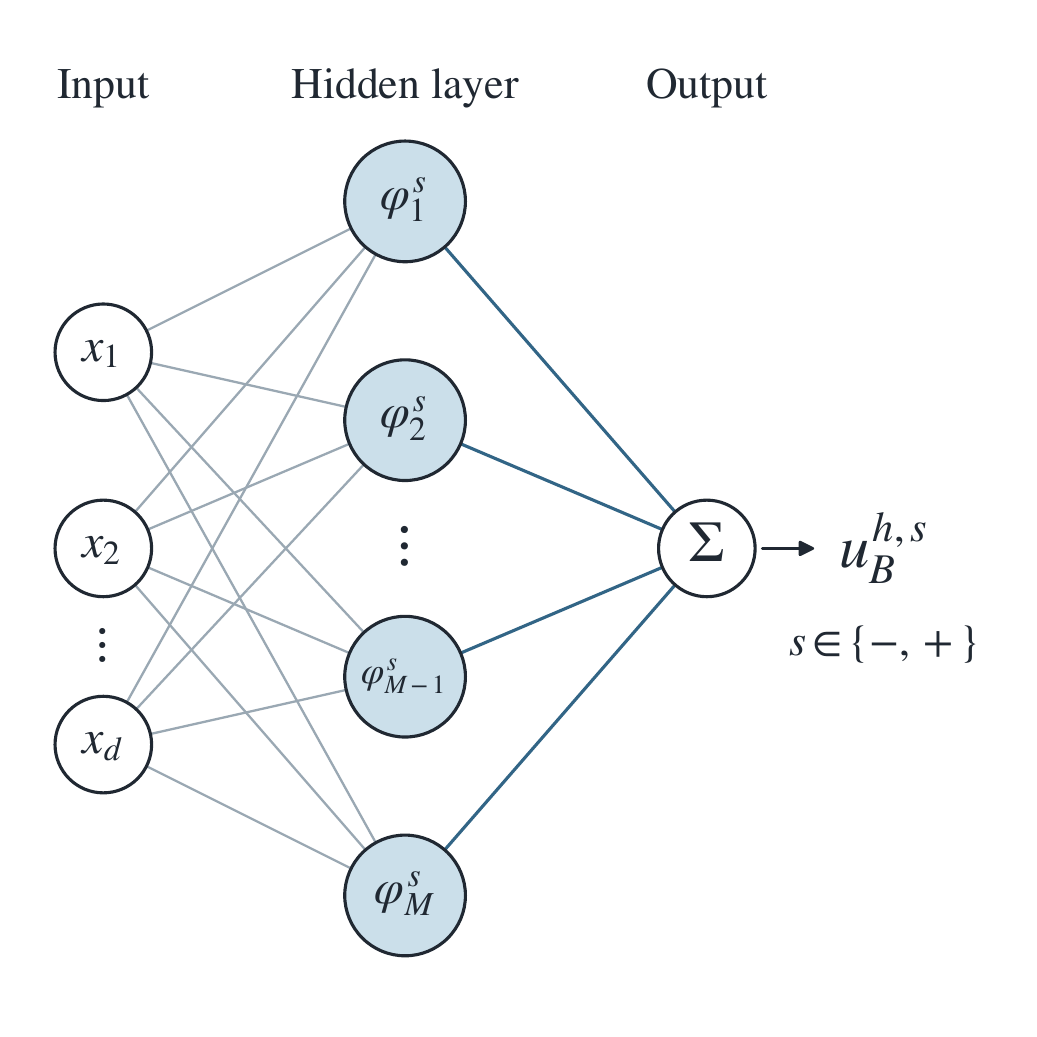}
        \caption{Local ELM architecture.}
        \label{fig:local_elm_architecture}
    \end{subfigure}
    \caption{Hybrid domain decomposition and local ELM architecture.}
    \label{fig:hybrid_framework}
\end{figure}

\subsection{Uniform Cartesian mesh and hybrid domain decomposition}
\label{subsec:mesh-band}

For the core formulation, we take the computational domain to be a
Cartesian box
\[
    \Omega=\prod_{r=1}^{d}(a_r,b_r)
\]
and introduce a uniform Cartesian background mesh with cell side length
$h$. Assuming that $N_r=(b_r-a_r)/h\in\mathbb N$ for each coordinate
direction, the mesh is given by
\begin{equation}
\label{eq:uniform-cartesian-mesh}
\begin{aligned}
    \mathcal T_h
    &=
    \{K_{\boldsymbol i}:1\leq i_r\leq N_r,
        \quad r=1,\ldots,d\},\\
    K_{\boldsymbol i}
    &=
    \prod_{r=1}^{d}
    \bigl(a_r+(i_r-1)h,\,a_r+i_rh\bigr),
\end{aligned}
\end{equation}
where $\boldsymbol i=(i_1,\ldots,i_d)$ is the cell index. The background
mesh is constructed independently of the physical interface $\Gamma$.
The cells intersected by the physical interface are collected into the
interface band:
\begin{equation}
\label{eq:band-cell-set}
    \mathcal T_h^{B}
    =
    \left\{
        K\in\mathcal T_h:
        \overline K\cap\Gamma\neq\emptyset
    \right\}.
\end{equation}
The use of the cell closure includes cells for which $\Gamma$ meets a
cell face, edge, or vertex. The remaining cells form the FEM-cell set
\begin{equation}
\label{eq:fem-cell-set}
    \mathcal T_h^{F}
    =
    \mathcal T_h\setminus\mathcal T_h^{B}.
\end{equation}
The corresponding open regions are
\begin{equation}
\label{eq:fem-band-domains}
\begin{aligned}
    \Omega_h^{B}
    =
    \operatorname{int}\!\left(
        \bigcup_{K\in\mathcal T_h^{B}}\overline K
    \right),\quad 
    \Omega_h^{F}
    =
    \operatorname{int}\!\left(
        \bigcup_{K\in\mathcal T_h^{F}}\overline K
    \right).
\end{aligned}
\end{equation}
Since every cell intersected by $\Gamma$ belongs to
$\mathcal T_h^{B}$, each cell in $\mathcal T_h^{F}$ lies entirely in
one of the physical subdomains. Hence, the FEM region contains only
uncut cells and can be treated by a standard conforming finite element
discretization.
The artificial coupling interface is defined by
\begin{equation}
\label{eq:artificial-interface}
    \Gamma_c
    =
    \partial\Omega_h^{F}
    \cap
    \partial\Omega_h^{B}
    \cap\Omega,
    \qquad
    \Gamma_c^s
    =
    \Gamma_c\cap\Omega^s,
    \qquad s\in\{-,+\}.
\end{equation}
The unit normal $\boldsymbol n_F$ on $\Gamma_c$ points from
$\Omega_h^{F}$ into $\Omega_h^{B}$. We further write
\begin{equation}
\label{eq:band-side-regions}
    \Omega_h^{B,s}
    =
    \Omega_h^{B}\cap\Omega^s,
    \qquad s\in\{-,+\}.
\end{equation}
The FEM–RFM partition is therefore obtained by classifying cells on a prescribed background mesh according to their intersection with the physical interface. This construction separates background mesh generation from the allocation of the two approximation methods and retains standard cell geometry and polynomial spaces throughout the FEM region. The numerical partition into \(\Omega_h^F\) and \(\Omega_h^B\) is distinct from the physical decomposition into \(\Omega^-\) and \(\Omega^+\): the prescribed jumps are imposed on \(\Gamma\), while the FEM and random-feature representations are coupled on \(\Gamma_c\). Changes in the physical interface affect the band selection, phase classification, and geometric data required for sampling and flux evaluation.
A representative hybrid decomposition is illustrated in
Figure~\ref{fig:hybrid_domain_decomposition}.

\begin{remark}[Enlargement of the interface band]
\label{rem:band-enlargement}
The definition~\eqref{eq:band-cell-set} gives the minimal interface
band. If desired, $\mathcal T_h^{B}$ may be enlarged locally by
including selected neighboring cells, for example when the interface
passes very close to a grid line or when a wider region is desirable
for the localized random-feature approximation. The subsequent
decomposition and coupling formulation remain unchanged after replacing
$\mathcal T_h^{B}$ by the enlarged cell set.
\end{remark}

If the interface band reaches the external boundary on a coarse mesh,
the same decomposition is retained. On the portion of $\partial\Omega$
covered by the FEM region, the Dirichlet data are imposed through the
standard finite element boundary degrees of freedom. On the portion
covered by the interface band, the condition $u=g$ is imposed directly
through boundary collocation residuals of the local random-feature
approximation. The band-interior PDE residuals and the coupling across
$\Gamma_c$ are otherwise unchanged.

\begin{remark}[Irregular outer boundaries]
\label{rem:irregular-boundary}
The Cartesian-box assumption above is used to describe the core
formulation. For a domain with an irregular outer boundary, $\Omega$
may be embedded in a larger Cartesian box and discretized by a uniform
background mesh. In addition to the interface band $\Omega_h^B$ associated with $\Gamma$, the cells intersected by $\partial\Omega$, optionally together with selected neighboring cells, are collected into a boundary band. Localized random-feature approximations are used on the portions of these bands inside $\Omega$, while standard finite elements are retained on cells lying entirely inside $\Omega$ and outside both bands.

The prescribed Dirichlet condition is imposed directly on the physical
boundary $\partial\Omega$, and the boundary-band random-feature
approximation is coupled to the interior FEM region across the resulting
grid-aligned artificial boundary using the same trace and variational
normal-flux coupling mechanism. The feasibility of this extension is
demonstrated numerically in
Section~\ref{subsec:irregular_outer_boundary}.
\end{remark}

\subsection{Discrete FEM and local ELM approximation spaces}
\label{subsec:discrete-spaces}

Let $p\geq1$. For each $K\in\mathcal T_h^F$, let $Q_p(K)$ denote the
standard tensor-product polynomial space of degree at most $p$ in each
coordinate. The conforming finite element space on the FEM region is
\begin{equation}
\label{eq:fem-space}
    V_h^{F}
    =
    \left\{
        v_h\in C^0(\overline{\Omega_h^{F}}):
        v_h|_K\in Q_p(K),
        \quad K\in\mathcal T_h^{F}
    \right\}.
\end{equation}

Let $\mathcal N_F=\{1,\ldots,n_F\}$ index the global Lagrange nodes
$\{\boldsymbol x_a\}_{a\in\mathcal N_F}$ and the corresponding nodal
basis functions $\{N_a\}_{a\in\mathcal N_F}$, satisfying
\[
    N_a(\boldsymbol x_b)=\delta_{ab}.
\]
The finite element approximation is written as
\begin{equation}
\label{eq:fem-expansion}
    u_F^h(\boldsymbol x)
    =
    \sum_{a\in\mathcal N_F}
    U_aN_a(\boldsymbol x),
    \qquad
    \boldsymbol U
    =
    (U_1,\ldots,U_{n_F})^T.
\end{equation}

The nodal indices are partitioned into the disjoint sets
\begin{equation}
\label{eq:fem-node-partition}
\begin{aligned}
    \mathcal D
    &=
    \{a\in\mathcal N_F:
      \boldsymbol x_a\in\partial\Omega\},\\
    \mathcal C
    &=
    \{a\in\mathcal N_F\setminus\mathcal D:
      \boldsymbol x_a\in\Gamma_c\},\\
    \mathcal I
    &=
    \mathcal N_F\setminus
    (\mathcal C\cup\mathcal D).
\end{aligned}
\end{equation}
Thus, boundary nodes take precedence over coupling nodes when
$\Gamma_c$ meets $\partial\Omega$. For
\[
    (X,\mathcal X)
    \in
    \{(I,\mathcal I),(C,\mathcal C),(D,\mathcal D)\},
\]
let $n_X=|\mathcal X|$ and
$\boldsymbol U_X=(U_a)_{a\in\mathcal X}$, using a fixed ordering within
each index set. The finite element Dirichlet values are imposed nodally:
\begin{equation}
\label{eq:fem-dirichlet-data}
    \boldsymbol U_D
    =
    \boldsymbol g_D,
    \qquad
    \boldsymbol g_D
    =
    \bigl(g(\boldsymbol x_a)\bigr)_{a\in\mathcal D}.
\end{equation}
The corresponding trial and test spaces are
\begin{equation}
\label{eq:fem-trial-test-spaces}
\begin{aligned}
    V_{h,g}^{F}
    &=
    \left\{
        v_h\in V_h^{F}:
        v_h(\boldsymbol x_a)=g(\boldsymbol x_a),
        \quad a\in\mathcal D
    \right\},\\
    V_{h,0}^{F}
    &=
    \left\{
        v_h\in V_h^{F}:
        v_h(\boldsymbol x_a)=0,
        \quad a\in\mathcal D
    \right\}.
\end{aligned}
\end{equation}
The coupling values $\boldsymbol U_C$ remain unknown.
For each physical side $s\in\{-,+\}$, we introduce $M$ fixed random
features
\begin{equation}
\label{eq:elm-feature}
    \varphi_j^s(\boldsymbol x)
    =
    \sigma\!\left(
        \boldsymbol\omega_j^s\cdot\boldsymbol x+b_j^s
    \right),
    \qquad
    j=1,\ldots,M,
\end{equation}
where the input weights
$\boldsymbol\omega_j^s\in\mathbb R^d$ and biases
$b_j^s\in\mathbb R$ are randomly initialized once and then frozen.
The activation function
$\sigma:\mathbb R\rightarrow\mathbb R$ is assumed to satisfy
$\sigma\in C^2(\mathbb R)$ so that the second-order differential
operator can be evaluated directly on the trial functions. Only the
output coefficients are treated as unknowns in the local ELM
representations.
The local ELM trial spaces are defined by
\begin{equation}
\label{eq:one-sided-elm-spaces}
    \mathcal E_M^s
    =
    \operatorname{span}
    \left\{
        \varphi_j^s|_{\Omega_h^{B,s}}:
        j=1,\ldots,M
    \right\},
    \qquad
    s\in\{-,+\}.
\end{equation}
The corresponding band approximations are
\begin{equation}
\label{eq:two-sided-elm}
    u_B^{h,s}(\boldsymbol x)
    =
    \sum_{j=1}^{M}
    c_j^s\varphi_j^s(\boldsymbol x),
    \qquad
    \boldsymbol x\in\Omega_h^{B,s},
    \qquad
    s\in\{-,+\}.
\end{equation}
Let
\begin{equation}
\label{eq:one-sided-elm-coefficients}
    \boldsymbol c^s
    =
    \begin{bmatrix}
        c_1^s&\cdots&c_M^s
    \end{bmatrix}^{T}
    \in\mathbb R^M,
    \qquad
    s\in\{-,+\},
\end{equation}
and collect the two coefficient vectors as
\begin{equation}
\label{eq:elm-coefficient-vector}
    \boldsymbol c
    =
    \begin{bmatrix}
        \boldsymbol c^+\\
        \boldsymbol c^-
    \end{bmatrix}
    \in\mathbb R^{2M}.
\end{equation}

For the subsequent coupled assembly, the two local feature vectors are
embedded into the common coefficient space $\mathbb R^{2M}$. Define
$\Phi^\pm(\boldsymbol x)\in\mathbb R^{1\times2M}$ by
\begin{equation}
\label{eq:embedded-feature-rows}
\begin{aligned}
    \Phi^+(\boldsymbol x)
    &=
    \begin{bmatrix}
        \varphi_1^+(\boldsymbol x)
        &
        \cdots
        &
        \varphi_M^+(\boldsymbol x)
        &
        \boldsymbol 0_{1\times M}
    \end{bmatrix},
    \\
    \Phi^-(\boldsymbol x)
    &=
    \begin{bmatrix}
        \boldsymbol 0_{1\times M}
        &
        \varphi_1^-(\boldsymbol x)
        &
        \cdots
        &
        \varphi_M^-(\boldsymbol x)
    \end{bmatrix}.
\end{aligned}
\end{equation}
The two band fields can then be written uniformly as
\begin{equation}
\label{eq:elm-compact-representation}
    u_B^{h,s}(\boldsymbol x)
    =
    \Phi^s(\boldsymbol x)\boldsymbol c,
    \qquad
    \boldsymbol x\in\Omega_h^{B,s},
    \qquad
    s\in\{-,+\}.
\end{equation}
The two local random-feature representations occupy disjoint
coefficient blocks. They are coupled through the prescribed jump
conditions on $\Gamma$, and each is coupled to the FEM field
through the trace and variational normal-flux residuals on
$\Gamma_c^s$, $s\in\{-,+\}$.

\subsection{FEM weak formulation and matrix assembly}
\label{subsec:fem-assembly}

On the FEM region $\Omega_h^F$, define the bilinear form and volume load functional
\begin{equation}
\label{eq:fem-bilinear-load}
    a_F(w_h,v_h)
    =
    \int_{\Omega_h^{F}}
    \left(
        \beta\nabla w_h\cdot\nabla v_h
        +
        \kappa w_hv_h
    \right)
    \,\mathrm d\boldsymbol x,
    \qquad
    \ell_F(v_h)
    =
    \int_{\Omega_h^{F}}
    f v_h\,\mathrm d\boldsymbol x.
\end{equation}
Integration by parts over $\Omega_h^F$ produces a boundary term on
$\Gamma_c$. Replacing the FEM-side normal flux there by the corresponding
band-side flux leads to the residual
\begin{equation}
\label{eq:fem-weak-residual}
\begin{aligned}
    \mathcal R_F(v_h;\boldsymbol U,\boldsymbol c)
    ={}&a_F(u_F^h,v_h)-\ell_F(v_h)
       -\sum_{s\in\{-,+\}}\int_{\Gamma_c^s}
          \beta^s\partial_{\boldsymbol n_F}u_B^{h,s}
          \,v_h\,\mathrm ds,
    \qquad v_h\in V_{h,0}^{F}.
\end{aligned}
\end{equation}
The sign of the boundary term follows from the orientation of
$\boldsymbol n_F$ as the outward normal of the FEM region.
For an interior basis function $N_a$, $a\in\mathcal I$, the trace of
$N_a$ vanishes on $\Gamma_c$. Hence,
\[
    \mathcal R_F(N_a;\boldsymbol U,\boldsymbol c)=0
\]
reduces to the standard interior Galerkin equation. These equations are
enforced algebraically and subsequently eliminated by the Schur
reduction. The basis functions associated with
$\mathcal C$ have nonzero traces on $\Gamma_c$; the corresponding
Galerkin rows contain the band-side flux contribution and are retained
as variational FEM--RFM coupling residuals.
The FEM matrix $K\in\mathbb R^{n_F\times n_F}$ and load vector
$\boldsymbol F\in\mathbb R^{n_F}$ are assembled cellwise as
\begin{equation}
\label{eq:fem-matrix-assembly}
\begin{aligned}
    K_{ab}
    &=
    \sum_{K\in\mathcal T_h^{F}}
    \int_K
    \left(
        \beta(\boldsymbol x)
        \nabla N_b\cdot\nabla N_a
        +
        \kappa N_bN_a
    \right)
    \,\mathrm d\boldsymbol x,
    \\
    F_a
    &=
    \sum_{K\in\mathcal T_h^{F}}
    \int_K
    f(\boldsymbol x)N_a(\boldsymbol x)
    \,\mathrm d\boldsymbol x,
    \qquad
    a,b\in\mathcal N_F.
\end{aligned}
\end{equation}
These are standard finite element volume contributions on uncut cells.
The zeroth-order term contributes only to the volume operator and therefore
does not modify the normal-flux coupling on $\Gamma_c$. The coupling flux
is assembled separately as a boundary contribution.
For $X,Y\in\{I,C,D\}$, use the block notation
\begin{equation}
\label{eq:fem-block-notation}
    K_{XY}=(K_{ab})_{a\in\mathcal X,\,b\in\mathcal Y}
        \in\mathbb R^{n_X\times n_Y},
    \qquad
    \boldsymbol F_X=(F_a)_{a\in\mathcal X}\in\mathbb R^{n_X}.
\end{equation}
For an interior basis function $N_a$, $a\in\mathcal I$, the trace of
$N_a$ vanishes on $\Gamma_c$. Therefore, the coupling flux contribution
in \eqref{eq:fem-weak-residual} is zero, and
$\mathcal R_F(N_a;\boldsymbol U,\boldsymbol c)=0$ reduces to the
interior Galerkin equations.
\begin{equation}
\label{eq:interior-fem-equation}
    K_{II}\boldsymbol U_I
    +K_{IC}\boldsymbol U_C
    +K_{ID}\boldsymbol g_D
    =\boldsymbol F_I.
\end{equation}
Once $\boldsymbol U_C$ and $\boldsymbol g_D$ are specified, this is a
discrete Dirichlet problem on each FEM component. With the coupling and
Dirichlet values treated as prescribed boundary data for the interior FEM
solve, the functions associated with $\mathcal I$ have zero trace on the
boundary of each FEM component. The uniform positivity of $\beta$ and the
assumption $\kappa\geq0$ therefore imply that $K_{II}$ is symmetric
positive definite.

\subsection{Schur reduction and FEM--RFM coupling on
\texorpdfstring{$\Gamma_c$}{Gamma-c}}
\label{subsec:fem-schur}

The coupling retains the FEM nodal values and the random-feature
coefficients as independent unknowns. Write the ordered coupling index
set as
\[
    \mathcal C=(c_1,\ldots,c_{n_C}),
\]
and let
\[
    \boldsymbol x_i^C=\boldsymbol x_{c_i},
    \qquad i=1,\ldots,n_C,
\]
denote the corresponding FEM coupling nodes on $\Gamma_c$. For each
$i$, let $s_i\in\{-,+\}$ denote the physical-side label of
$\boldsymbol x_i^C$, so that
\[
    \boldsymbol x_i^C\in\Omega^{s_i}.
\]
Define the trace matrix
$\Phi_C\in\mathbb R^{n_C\times2M}$ by
\begin{equation}
\label{eq:coupling-trace-matrix}
    (\Phi_C)_{ij}
    =\bigl[\Phi^{s_i}(\boldsymbol x_i^C)\bigr]_j,
    \qquad i=1,\ldots,n_C,\quad j=1,\ldots,2M.
\end{equation}
Here $[\Phi^s]_j$ denotes the $j$th component of the embedded feature row.
The nodal trace residual is therefore
\begin{equation}
\label{eq:trace-residual}
    \boldsymbol r_{\mathrm{tr}}
    =\boldsymbol U_C-\Phi_C\boldsymbol c.
\end{equation}
Thus, $\boldsymbol r_{\mathrm{tr}}$ measures the nodal mismatch between
the FEM trace and the local random-feature trace on $\Gamma_c$. This
continuity condition is imposed in the least-squares sense.

For the coupling test function $N_{c_i}$, the band flux contribution in
\eqref{eq:fem-weak-residual} is linear in $\boldsymbol c$. Define
$A_C\in\mathbb R^{n_C\times2M}$ by
\begin{equation}
\label{eq:band-flux-load}
    (A_C)_{ij}
    =\sum_{s\in\{-,+\}}\int_{\Gamma_c^s}
       \beta^s(\boldsymbol x)
       \partial_{\boldsymbol n_F}[\Phi^s(\boldsymbol x)]_j
       \,N_{c_i}(\boldsymbol x)\,\mathrm ds.
\end{equation}
Thus, $(A_C\boldsymbol c)_i$ is the band-side normal flux integrated
against the FEM trace basis $N_{c_i}|_{\Gamma_c}$. These integrals are
assembled over the grid-aligned faces of $\Gamma_c$, with
$\partial_{\boldsymbol n_F}q=\nabla q\cdot\boldsymbol n_F$.
Before eliminating the interior unknowns, the coupling residual is
\begin{equation}
\label{eq:uncondensed-flux-residual}
    \boldsymbol r_c
    =K_{CI}\boldsymbol U_I+K_{CC}\boldsymbol U_C
      +K_{CD}\boldsymbol g_D-A_C\boldsymbol c-\boldsymbol F_C.
\end{equation}

To eliminate the interior FEM degrees of freedom, define
$\boldsymbol U_I^0\in\mathbb R^{n_I}$ and
$W\in\mathbb R^{n_I\times n_C}$ as the solutions of
\begin{equation}
\label{eq:interior-response-problems}
    K_{II}\boldsymbol U_I^0
    =\boldsymbol F_I-K_{ID}\boldsymbol g_D,
    \qquad
    K_{II}W=K_{IC}.
\end{equation}
Equation~\eqref{eq:interior-fem-equation} then gives
\begin{equation}
\label{eq:fem-interior-elimination}
    \boldsymbol U_I
    =\boldsymbol U_I^0-W\boldsymbol U_C.
\end{equation}
Substitution into \eqref{eq:uncondensed-flux-residual} yields the Schur
complement and condensed load
\begin{equation}
\label{eq:schur-quantities}
\begin{aligned}
    S_{CC}
    &=K_{CC}-K_{CI}W
        \in\mathbb R^{n_C\times n_C},\\
    \boldsymbol b_C
    &=\boldsymbol F_C-K_{CI}\boldsymbol U_I^0
                      -K_{CD}\boldsymbol g_D
        \in\mathbb R^{n_C}.
\end{aligned}
\end{equation}
Hence,
\begin{equation}
\label{eq:flux-coupling-residual}
    \boldsymbol r_c
    =S_{CC}\boldsymbol U_C-A_C\boldsymbol c-\boldsymbol b_C.
\end{equation}
The quantity $S_{CC}\boldsymbol U_C-\boldsymbol b_C$ is the condensed
FEM boundary response. The variational residual compares this response
with the band flux moments $A_C\boldsymbol c$. Together with the nodal trace values
$\Phi_C\boldsymbol c$, these flux moments provide the quantities
through which the band representation enters the FEM--RFM coupling.
For a fixed FEM region and linear volume operator, the interior FEM
discretization is represented in the reduced coupling by $S_{CC}$ and
$\boldsymbol b_C$. This separates the interior FEM solve from the boundary quantities
through which it communicates with the band. The interior equations remain satisfied through
\eqref{eq:fem-interior-elimination}, while both coupling residuals
$\boldsymbol r_{\mathrm{tr}}$ and $\boldsymbol r_c$ are minimized in
the reduced least-squares problem.

\subsection{Band residuals and physical interface conditions}
\label{subsec:residual-equations}

On each physical side, the band PDE residual is constructed from the
RFM approximation in $\Omega_h^{B,s}$. In the numerical
experiments, we use a strong-form collocation residual evaluated at
randomly sampled points. These points are generated independently
within each physical-side band region, and the resulting pointwise
residuals are incorporated into the reduced least-squares system.
A weak-form variant based on the same trial space is examined in
Appendix~\ref{app:weak_elm_fem}. Let
\begin{equation}
\label{eq:band-collocation-sets}
    \mathcal X_B^s
    =\{\boldsymbol x_m^s\}_{m=1}^{N_B^s}
      \subset\Omega_h^{B,s},
    \qquad s\in\{-,+\}.
\end{equation}
Define the differential operator
\begin{equation}
\label{eq:band-differential-operator}
    \mathcal L^s q
    =
    -\nabla\cdot(\beta^s\nabla q)
    +\kappa q,
    \qquad s\in\{-,+\}.
\end{equation}
The operator residual is evaluated at the collocation points.
Let
$A_{\mathrm{pde}}^s\in\mathbb R^{N_B^s\times2M}$ and
$\boldsymbol b_{\mathrm{pde}}^s\in\mathbb R^{N_B^s}$ be defined by
\begin{equation}
\label{eq:one-sided-pde-matrices}
    (A_{\mathrm{pde}}^s)_{mj}
    =\bigl(\mathcal L^s[\Phi^s]_j\bigr)(\boldsymbol x_m^s),
    \qquad
    (\boldsymbol b_{\mathrm{pde}}^s)_m
    =f^s(\boldsymbol x_m^s).
\end{equation}
Stacking the two sides gives
\begin{equation}
\label{eq:band-pde-matrix}
    A_{\mathrm{pde}}
    =\begin{bmatrix}A_{\mathrm{pde}}^+\\A_{\mathrm{pde}}^-\end{bmatrix},
    \qquad
    \boldsymbol b_{\mathrm{pde}}
    =\begin{bmatrix}
        \boldsymbol b_{\mathrm{pde}}^+\\
        \boldsymbol b_{\mathrm{pde}}^-
      \end{bmatrix},
    \qquad N_B=N_B^++N_B^-.
\end{equation}
Thus $A_{\mathrm{pde}}\in\mathbb R^{N_B\times2M}$, and the band PDE
residual is
\begin{equation}
\label{eq:band-pde-residual}
    \boldsymbol r_{\mathrm{pde}}
    =A_{\mathrm{pde}}\boldsymbol c-\boldsymbol b_{\mathrm{pde}}.
\end{equation}

The physical jump conditions are imposed independently of the coupling
conditions on $\Gamma_c$. Choose interface points
$\mathcal X_\Gamma=\{\boldsymbol x_\ell^\Gamma\}_{\ell=1}^{N_\Gamma}$
at which the unit normal is defined, and write
$\boldsymbol n_\ell=\boldsymbol n(\boldsymbol x_\ell^\Gamma)$.
Define the value-jump matrix $A_w\in\mathbb R^{N_\Gamma\times2M}$ and
its data vector by
\begin{equation}
\label{eq:value-jump-matrix}
\begin{aligned}
    (A_w)_{\ell j}
    &=[\Phi^+(\boldsymbol x_\ell^\Gamma)]_j
      -[\Phi^-(\boldsymbol x_\ell^\Gamma)]_j,\\
    (\boldsymbol b_w)_\ell
    &=w(\boldsymbol x_\ell^\Gamma).
\end{aligned}
\end{equation}
Similarly, the flux-jump matrix $A_v\in\mathbb R^{N_\Gamma\times2M}$
and its data vector are
\begin{equation}
\label{eq:flux-jump-matrix}
\begin{aligned}
    (A_v)_{\ell j}
    ={}&\beta^+(\boldsymbol x_\ell^\Gamma)
        \nabla[\Phi^+]_j(\boldsymbol x_\ell^\Gamma)
        \cdot\boldsymbol n_\ell\\
       &-\beta^-(\boldsymbol x_\ell^\Gamma)
        \nabla[\Phi^-]_j(\boldsymbol x_\ell^\Gamma)
        \cdot\boldsymbol n_\ell,\\
    (\boldsymbol b_v)_\ell
    ={}&v(\boldsymbol x_\ell^\Gamma).
\end{aligned}
\end{equation}
Both terms use the same physical normal, oriented from $\Omega^-$ to
$\Omega^+$. The two jump residuals are combined as
\begin{equation}
\label{eq:physical-interface-matrix}
    A_\Gamma
    =\begin{bmatrix}A_w\\A_v\end{bmatrix}
      \in\mathbb R^{2N_\Gamma\times2M},
    \qquad
    \boldsymbol b_\Gamma
    =\begin{bmatrix}\boldsymbol b_w\\\boldsymbol b_v\end{bmatrix},
    \qquad
    \boldsymbol r_\Gamma
    =A_\Gamma\boldsymbol c-\boldsymbol b_\Gamma.
\end{equation}
The physical interface therefore couples the two local random-feature
representations through the prescribed solution and flux jumps.
To include the case in which the interface band reaches the external
boundary, let
\[
    \mathcal X_{\mathrm{bd}}^s
    =
    \{\boldsymbol x_r^{\mathrm{bd},s}\}_{r=1}^{N_{\mathrm{bd}}^s}
    \subset
    \partial\Omega\cap\overline{\Omega_h^{B,s}},
    \qquad s\in\{-,+\},
\]
denote the corresponding boundary collocation points. Define
$A_{\mathrm{bd}}^s\in\mathbb R^{N_{\mathrm{bd}}^s\times2M}$ and
$\boldsymbol b_{\mathrm{bd}}^s\in\mathbb R^{N_{\mathrm{bd}}^s}$ by
\begin{equation}
\label{eq:boundary-residual-matrices}
    (A_{\mathrm{bd}}^s)_{rj}
    =
    [\Phi^s(\boldsymbol x_r^{\mathrm{bd},s})]_j,
    \qquad
    (\boldsymbol b_{\mathrm{bd}}^s)_r
    =
    g(\boldsymbol x_r^{\mathrm{bd},s}).
\end{equation}
Stacking the two physical sides gives
\begin{equation}
\label{eq:boundary-residual}
    A_{\mathrm{bd}}
    =
    \begin{bmatrix}
        A_{\mathrm{bd}}^+\\
        A_{\mathrm{bd}}^-
    \end{bmatrix},
    \qquad
    \boldsymbol b_{\mathrm{bd}}
    =
    \begin{bmatrix}
        \boldsymbol b_{\mathrm{bd}}^+\\
        \boldsymbol b_{\mathrm{bd}}^-
    \end{bmatrix},
    \qquad
    \boldsymbol r_{\mathrm{bd}}
    =
    A_{\mathrm{bd}}\boldsymbol c-\boldsymbol b_{\mathrm{bd}},
\end{equation}
where
\[
    N_{\mathrm{bd}}
    =
    N_{\mathrm{bd}}^+
    +
    N_{\mathrm{bd}}^-.
\]
If the interface band does not intersect $\partial\Omega$, then
$N_{\mathrm{bd}}=0$, and the corresponding boundary-residual block is
absent.

\subsection{Reduced least-squares formulation and reconstruction}
\label{subsec:augmented-ls-system}

After imposing the Dirichlet data and eliminating the interior FEM
unknowns, the reduced coefficient vector is
\begin{equation}
\label{eq:reduced-unknown}
    \boldsymbol z
    =
    \begin{bmatrix}
        \boldsymbol U_C\\
        \boldsymbol c
    \end{bmatrix}
    \in\mathbb R^{n_C+2M}.
\end{equation}

The trace, variational flux, band PDE, physical interface, and, when
present, band-boundary residuals form the block system
\begin{equation}
\label{eq:augmented-matrix}
    A
    =
    \begin{bmatrix}
        I_C&-\Phi_C\\
        S_{CC}&-A_C\\
        0&A_{\mathrm{pde}}\\
        0&A_\Gamma\\
        0&A_{\mathrm{bd}}
    \end{bmatrix},
    \qquad
    \boldsymbol b
    =
    \begin{bmatrix}
        \boldsymbol0\\
        \boldsymbol b_C\\
        \boldsymbol b_{\mathrm{pde}}\\
        \boldsymbol b_\Gamma\\
        \boldsymbol b_{\mathrm{bd}}
    \end{bmatrix}.
\end{equation}
Here $I_C$ is the $n_C\times n_C$ identity matrix and the zero blocks
have the corresponding dimensions. In particular,
\[
    A\in\mathbb R^{n_R\times(n_C+2M)},
    \qquad
    n_R
    =
    2n_C+N_B+2N_\Gamma+N_{\mathrm{bd}},
\]
and
\begin{equation}
\label{eq:assembled-residual-vector}
    A\boldsymbol z-\boldsymbol b
    =
    \begin{bmatrix}
        \boldsymbol r_{\mathrm{tr}}\\
        \boldsymbol r_c\\
        \boldsymbol r_{\mathrm{pde}}\\
        \boldsymbol r_\Gamma\\
        \boldsymbol r_{\mathrm{bd}}
    \end{bmatrix}.
\end{equation}

We use row equilibration to balance the algebraic scales of the residual
equations. Denoting the $i$th row of $A$ by $A_{i,:}$, define
\begin{equation}
\label{eq:row-equilibration}
    D_r
    =
    \operatorname{diag}(d_1,\ldots,d_{n_R}),
    \qquad
    d_i
    =
    \begin{cases}
        \|A_{i,:}\|_2^{-1},
        & \|A_{i,:}\|_2>0,\\
        1,
        & \|A_{i,:}\|_2=0.
    \end{cases}
\end{equation}
The reduced coefficients are determined by
\begin{equation}
\label{eq:augmented-ls-problem}
    \boldsymbol z_h
    \in
    \operatorname*{arg\,min}_{\boldsymbol z\in\mathbb R^{n_C+2M}}
    \|D_r(A\boldsymbol z-\boldsymbol b)\|_2^2.
\end{equation}
If the minimizer is not unique, $\boldsymbol z_h$ denotes the minimizer
with the smallest Euclidean norm. The interior FEM equations are incorporated as equality
constraints through the Schur reduction, and
\eqref{eq:augmented-ls-problem} determines the remaining $n_C+2M$
unknowns.

Writing
$\boldsymbol z_h=((\boldsymbol U_C^h)^T,(\boldsymbol c^h)^T)^T$,
recover the remaining FEM coefficients from
\begin{equation}
\label{eq:hybrid-coefficient-reconstruction}
    \boldsymbol U_I^h
    =
    \boldsymbol U_I^0-W\boldsymbol U_C^h,
    \qquad
    \boldsymbol U_D^h
    =
    \boldsymbol g_D.
\end{equation}
The resulting composite solution is defined almost everywhere in
$\Omega$ by
\begin{equation}
\label{eq:hybrid-solution-reconstruction}
    u_h(\boldsymbol x)
    =
    \begin{cases}
        \displaystyle
        \sum_{a\in\mathcal N_F}
        U_a^hN_a(\boldsymbol x),
        & \boldsymbol x\in\Omega_h^{F},\\[2mm]
        \Phi^s(\boldsymbol x)\boldsymbol c^h,
        & \boldsymbol x\in\Omega_h^{B,s},
        \quad s\in\{-,+\}.
    \end{cases}
\end{equation}

In the present formulation, trace continuity on $\Gamma_c$ is enforced
through the discrete least-squares residual
$\boldsymbol r_{\mathrm{tr}}$. Therefore, the FEM and random-feature
traces are not required to coincide exactly on $\Gamma_c$, and the
resulting composite approximation is generally nonconforming across the
artificial coupling interface. The proposed least-squares trace coupling
retains $\boldsymbol U_C$ and $\boldsymbol c$ as independent unknowns,
allowing the trace mismatch to be balanced together with the other
residuals in the reduced least-squares problem. Appendix~
\ref{app:star_supplementary} compares the proposed soft trace coupling
with two hard alternatives based on nodal and moment equality
constraints.

\begin{remark}[Role of the Schur reduction]
\label{rem:schur-reduction}
The Schur reduction separates the size of the reconstructed
approximation from the number of coefficients entering the coupled
least-squares solve. After prescribing the Dirichlet values, the
uncondensed representation contains $n_I+n_C+2M$ coefficients.
The interior response problems in
\eqref{eq:interior-response-problems} are solved using a sparse
factorization of $K_{II}$ and the associated triangular solves,
leaving $n_C+2M$ unknowns in the reduced least-squares problem.
The interior coefficients are still present in the reconstructed FEM
field; they are determined by $\boldsymbol U_C$ through
\eqref{eq:fem-interior-elimination}.

For the assembled linear system, each coupling vector
$\boldsymbol U_C$ uniquely determines the interior FEM coefficients
through \eqref{eq:fem-interior-elimination}. The interior Galerkin
equations \eqref{eq:interior-fem-equation} are then satisfied in
exact arithmetic. In floating-point computation, the resulting
interior residual arises from numerical errors in the sparse
linear solves and coefficient reconstruction.
The reduced least-squares problem minimizes the remaining band
and interface residuals, including the trace and variational flux
residuals on $\Gamma_c$. Errors in $\boldsymbol U_C$ can still
propagate to the reconstructed FEM field through
\eqref{eq:fem-interior-elimination}.

Nominal coefficient counts alone do not characterize the computational
cost. The reduced matrix has $n_R$ rows, and the computation also
includes the formation of $W$ and $S_{CC}$, which are generally dense.
The reduction combines sparse interior solves with a dense reduced
least-squares problem. 
Its algebraic and numerical effects are
examined in Section~\ref{subsec:circular_high_contrast} by comparison
with a full least-squares formulation that includes the interior
Galerkin equations among the residual rows.
\end{remark}
\begin{remark}[Local subdivision of the random-feature band]
\label{rem:local-rfm-subdivision}
Each one-sided band may be further
partitioned into nonoverlapping subdomains according to its
connectivity, geometric complexity, or approximation requirements.
This local decomposition is consistent with domain-decomposition
strategies used in neural-network-based PDE solvers, where independent
local approximations are assigned to different subdomains and coupled
through interface conditions
\cite{shukla2021parallel}.
An independent set of fixed random features and output coefficients
is assigned to each subdomain. Across artificial boundaries between
adjacent subdomains within the same physical phase, continuity of
the solution and normal flux is imposed through additional
least-squares residuals. These conditions are distinct from the
prescribed jump conditions on the physical interface $\Gamma$.
The trace and variational normal-flux coupling with the FEM region
on $\Gamma_c$ is retained, and the local output coefficients are
collected into a common coefficient vector. For a fixed FEM region
and linear volume operator, this subdivision changes the band
representation and introduces additional residuals in the band, but
leaves the interior FEM matrix unchanged. The same Schur reduction
therefore applies, with the revised band representation supplying
the required trace values and flux moments. This further subdivision
is used for the bicontinuous composite in
Section~\ref{subsec:bicontinuous_homogenization}.
\end{remark}
\subsection{Extensions of the hybrid formulation}
\label{subsec:method_extensions}

The spatial decomposition and FEM--RFM coupling can also be used
for the following extensions.

\paragraph{Quasilinear elliptic interface problems.}
For quasilinear elliptic problems with a solution-dependent diffusion
coefficient $\beta^s(\boldsymbol x,u^s)$, the same spatial decomposition
and coupling structure are retained. The FEM volume operator, band
PDE residuals, and flux terms depend on the approximate solution,
while the trace residual retains its linear form. The resulting
coefficient problem is nonlinear. A numerical example is presented
in Section~\ref{subsec:quasilinear_circular_interface}.

\paragraph{Parabolic problems with prescribed moving interfaces.}
After implicit time discretization, a linear parabolic interface
problem leads at each time level to an elliptic problem of the form
\eqref{eq:linear-interface-problem}, with a zeroth-order term and a
source contribution from previous time levels. The same FEM--RFM
spatial discretization and coupling can then be applied. For a
prescribed moving interface $\Gamma(t)$, the geometric data are
evaluated at the current time level, and the history terms use
extensions of the corresponding phase-wise approximations where
needed. This extension is illustrated in
Section~\ref{subsec:moving_interface_parabolic}.

\section{Numerical examples}
\label{sec:numerical_examples}

The numerical experiments examine the mesh-refinement behavior of the
proposed hybrid FEM--RFM formulation, the algebraic role of the Schur reduction, and
the accuracy--cost tradeoff relative to full-subdomain RFM.
The additional geometric, quasilinear, and prescribed moving-interface
examples examine whether the same localized spatial decomposition and
coupling mechanism can be retained under changes in geometry and governing equation.

Throughout this section, absolute and relative errors are defined using a
common notation. For the solution error in the standard $L^2(\Omega)$ norm,
we set
\[
E_{L^2}^{\mathrm{abs}}
=
\|u-u_h\|_{L^2(\Omega)},
\qquad
E_{L^2}^{\mathrm{rel}}
=
\frac{\|u-u_h\|_{L^2(\Omega)}}
     {\|u\|_{L^2(\Omega)}}.
\]

Let $\mathcal P_h$ denote the partition consisting of the FEM cells and
the local random-feature subdomains, with each random-feature
subdomain restricted to its physical phase. This partition resolves
the physical interfaces, the FEM--RFM coupling interfaces, and any
artificial boundaries introduced by further subdivision of the band.
For a function $v$ that belongs to $H^1(D)$ on each $D\in\mathcal P_h$,
define
\[
|v|_{H^1_{\mathrm{br}}(\Omega)}
=
\left(
    \sum_{D\in\mathcal P_h}
    \|\nabla(v|_D)\|_{L^2(D)}^2
\right)^{1/2}.
\]
Here, the subscript $\mathrm{br}$ refers to the actual hybrid
partition. This seminorm measures the piecewise volume gradients;
it does not include trace-jump terms on the artificial interfaces.
The corresponding absolute and relative errors are
\[
E_{H^1}^{\mathrm{abs}}
=
|u-u_h|_{H^1_{\mathrm{br}}(\Omega)},
\qquad
E_{H^1}^{\mathrm{rel}}
=
\frac{|u-u_h|_{H^1_{\mathrm{br}}(\Omega)}}
     {|u|_{H^1_{\mathrm{br}}(\Omega)}}.
\]

For error measurements restricted to the interface band
$\Omega_h^B$, we further define
\[
E_{L^2,B}^{\mathrm{abs}}
=
\|u-u_h\|_{L^2(\Omega_h^B)},
\qquad
E_{H^1,B}^{\mathrm{abs}}
=
|u-u_h|_{H^1_{\mathrm{br}}(\Omega_h^B)}.
\]
The band seminorm uses the same partition restricted to $\Omega_h^B$,
including any artificial boundaries between local random-feature
subdomains within a physical phase.
When a nodal maximum error is reported, we use the discrete $L^\infty$
metric
\[
\|v\|_{L^\infty(\mathcal{N}_h)}
=
\max_{\boldsymbol{x}\in\mathcal{N}_h}
|v(\boldsymbol{x})|,
\]
where $\mathcal{N}_h$ is the evaluation-node set specified in the
corresponding experiment. We then define
\[
E_{L^\infty}^{\mathrm{abs}}
=
\|u-u_h\|_{L^\infty(\mathcal{N}_h)},
\qquad
E_{L^\infty}^{\mathrm{rel}}
=
\frac{\|u-u_h\|_{L^\infty(\mathcal{N}_h)}}
     {\|u\|_{L^\infty(\mathcal{N}_h)}}.
\]
If two numerical traces are evaluated at the same FEM--RFM coupling node,
they are treated as separate entries in $\mathcal{N}_h$. These definitions
are used throughout this section, independently of the notation adopted in
the cited references.

For two consecutive mesh sizes $h_{k-1}>h_k$, the observed convergence
order associated with an error measure $E(h)$ is defined by
\[
    \operatorname{Order}_k
    =
    \frac{
        \log\!\bigl(E(h_{k-1})/E(h_k)\bigr)
    }{
        \log\!\bigl(h_{k-1}/h_k\bigr)
    }.
\]
When the mesh size is successively halved, this reduces to
$\operatorname{Order}_k=\log_2(E(h_{k-1})/E(h_k))$.
The same definition is used for the error measures reported below. Here, $E(h)$ denotes any of the error measures defined above.

Unless otherwise specified, the zeroth-order coefficient is set to
$\kappa=0$ in the stationary elliptic examples. In all local ELM
realizations, each component of the random weight vector is sampled
independently from the uniform distribution
$\mathcal U[-2\pi,2\pi]$, while the bias is sampled from
$\mathcal U[-\pi,\pi]$. The hidden parameters are generated once and
then kept fixed during the corresponding computation. Except for the
feature-family comparison in
Section~\ref{subsec:circular_high_contrast}, the sinusoidal activation
$\sigma(\cdot)=\sin (\cdot)$ is used throughout the numerical experiments; the
hyperbolic-tangent activation is considered only in that comparison.

For the two-dimensional benchmark problems, each one-sided
random-feature representation uses $M=200$ fixed features and
$N_B^s=250$ band-interior PDE collocation points. These collocation
sets are intentionally chosen sufficiently large so that the localized
random-feature approximation does not become the dominant limitation
in the high-order mesh-refinement studies; the reported values are not
intended as minimal requirements of the method. In the three-dimensional
experiments, the corresponding numbers are increased to $M=1500$ and
$N_B^s=2000$ on each physical side. For the bicontinuous
homogenization problem, where the interface band is further decomposed
into multiple local random-feature subdomains, each local subnetwork
uses $M=50$ fixed features and $50$ band-interior PDE collocation
points. The above collocation counts refer only to the PDE residuals in the random-feature regions; the sampling of the physical interface $\Gamma$ and the quadrature used for the coupling interface $\Gamma_c$ are treated separately.

All numerical experiments were performed on a Linux compute node equipped with an Intel Xeon E5-2680 v4 processor (2.40 GHz, 28 logical cores) and 125 GB RAM. The implementations were executed with Python 3.12.3 and NumPy 1.26.4, linked against Intel MKL 2023.1 for the dense linear-algebra operations. No GPU acceleration was used.
\subsection{A high-contrast circular interface problem}
\label{subsec:circular_high_contrast}

We first consider a high-contrast circular interface problem studied in
\cite{ewing1999immersed,li2004immersed} and subsequently used in
interface optimal-control studies \cite{zhang2015immersed,lai2025hard}.
The main text uses this benchmark to examine mesh-refinement behavior under
large coefficient contrasts, the need for residual equilibration, and the
algebraic role of the Schur reduction. Additional parameter-sensitivity
results are reported in Appendix~\ref{app:circular_diagnostics}.

Let $\Omega=(-1,1)^2$ and $r_0=1/2$. The interface and the corresponding
subdomains are defined by
\[
\Gamma=
\left\{
    \boldsymbol{x}=(x_1,x_2)\in\Omega:
    x_1^2+x_2^2=r_0^2
\right\},\quad 
\Omega^-
=
\left\{
    \boldsymbol{x}\in\Omega:
    x_1^2+x_2^2<r_0^2
\right\},
\quad
\Omega^+
=
\left\{
    \boldsymbol{x}\in\Omega:
    x_1^2+x_2^2>r_0^2
\right\}.
\]
The interface problem is
\[
\begin{cases}
\displaystyle
-\nabla\cdot\bigl(\beta\nabla u\bigr)
=
-9\sqrt{x_1^2+x_2^2},
& \text{in }\Omega^-\cup\Omega^+,
\\[2mm]
\displaystyle
\jump{u}=0,
\qquad
\jump{\beta\partial_{\boldsymbol n}u}=0,
& \text{on }\Gamma,
\\[2mm]
u=g,
& \text{on }\partial\Omega,
\end{cases}
\]
where $\beta=\beta^\pm$ in $\Omega^\pm$. The exact solution is prescribed by
\[
    u(\boldsymbol{x})
    =
    \begin{cases}
        \displaystyle
        \frac{(x_1^2+x_2^2)^{3/2}}{\beta^-},
        & \boldsymbol{x}\in\Omega^-,
        \\[3mm]
        \displaystyle
        \frac{(x_1^2+x_2^2)^{3/2}}{\beta^+}
        +
        \left(
            \frac{1}{\beta^-}
            -
            \frac{1}{\beta^+}
        \right)r_0^3,
        & \boldsymbol{x}\in\Omega^+,
    \end{cases}
\]
with $g=u|_{\partial\Omega}$.
Unless otherwise stated, all computations in this example use the
row-equilibrated reduced least-squares formulation introduced in
Section~\ref{subsec:augmented-ls-system}, namely
\eqref{eq:row-equilibration}--\eqref{eq:augmented-ls-problem}.
The residual-scaling experiment below deliberately disables this
equilibration in order to isolate the algebraic effect of disparate
residual scales.

\begin{figure}[htbp]
    \centering
    \begin{subfigure}[t]{0.485\textwidth}
        \centering
        \includegraphics[width=\linewidth]{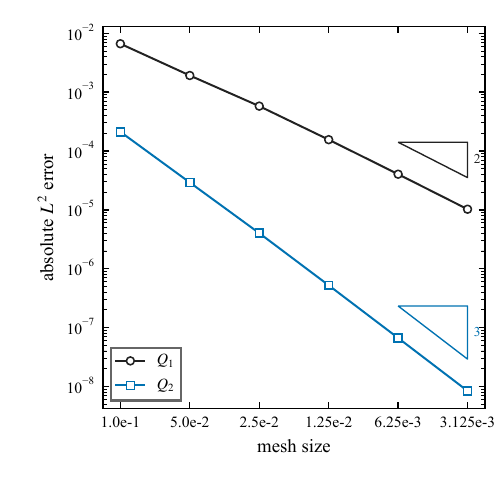}
        \caption{$E_{L^2}^{\mathrm{abs}}$, $(\beta^+,\beta^-)=(10^{12},1)$.}
        \label{fig:circular_high_contrast_l2_positive}
    \end{subfigure}
    \hfill
    \begin{subfigure}[t]{0.485\textwidth}
        \centering
        \includegraphics[width=\linewidth]{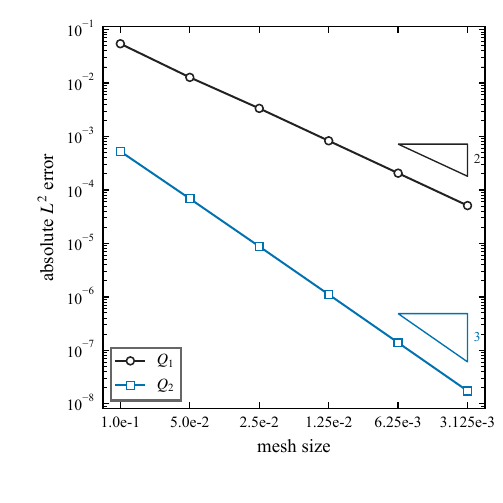}
        \caption{$E_{L^2}^{\mathrm{abs}}$, $(\beta^+,\beta^-)=(1,10^{12})$.}
        \label{fig:circular_high_contrast_l2_reversed}
    \end{subfigure}

    \begin{subfigure}[t]{0.485\textwidth}
        \centering
        \includegraphics[width=\linewidth]{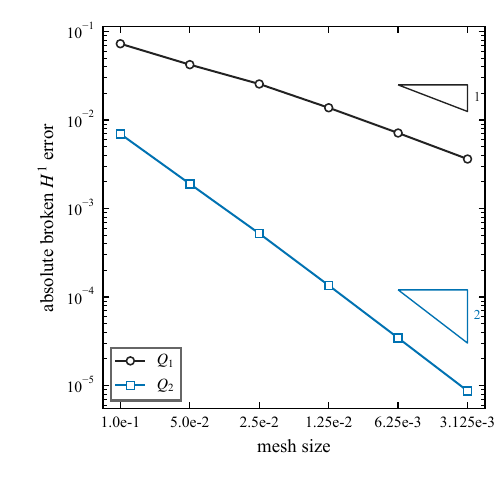}
        \caption{$E_{H^1}^{\mathrm{abs}}$, $(\beta^+,\beta^-)=(10^{12},1)$.}
        \label{fig:circular_high_contrast_h1_positive}
    \end{subfigure}
    \hfill
    \begin{subfigure}[t]{0.485\textwidth}
        \centering
        \includegraphics[width=\linewidth]{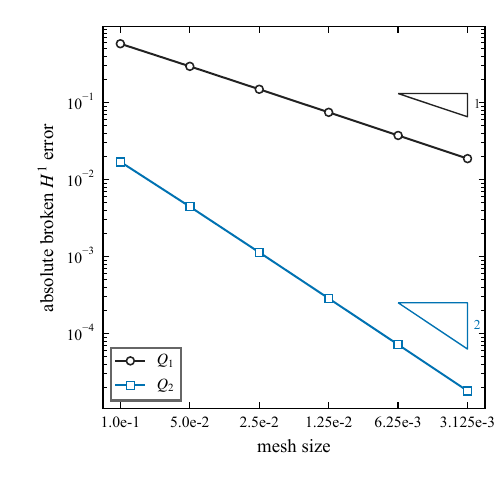}
        \caption{$E_{H^1}^{\mathrm{abs}}$, $(\beta^+,\beta^-)=(1,10^{12})$.}
        \label{fig:circular_high_contrast_h1_reversed}
    \end{subfigure}
    \caption{Mesh-refinement errors for the high-contrast circular interface problem.}
    \label{fig:circular_high_contrast_convergence}
\end{figure}

\paragraph{Mesh refinement under high coefficient contrasts.}
We consider $(\beta^+,\beta^-)=(10^{12},1)$ and
$(\beta^+,\beta^-)=(1,10^{12})$. The two configurations have the same
contrast magnitude but place the high-diffusivity phase on opposite sides
of the physical interface. Figure~\ref{fig:circular_high_contrast_convergence}
reports the $Q_1$ and $Q_2$ mesh-refinement results.
For both coefficient configurations, the $Q_1$ curves approach second-order
convergence in $E_{L^2}^{\mathrm{abs}}$ and first-order convergence in
$E_{H^1}^{\mathrm{abs}}$, while the $Q_2$ curves closely follow the
corresponding third- and second-order reference slopes.

\paragraph{Residual scaling and row equilibration.}
Before row equilibration, the flux-jump rows can be much larger than the
other residual equations because they contain the diffusion coefficients
explicitly. To isolate this effect, we disable row equilibration and multiply
only the physical flux-jump residual by a diagnostic factor
$s\in\{1,10^{-6},10^{-12},10^{-15}\}$ for the $Q_2$ computation with
$N=128$ and $h=1/64$. The complete data are given in
Table~\ref{tab:circular_high_contrast_scaling} of
Appendix~\ref{app:circular_diagnostics}.

For both coefficient orientations, reducing $s$ from $1$ to $10^{-12}$
reduces the spectral magnitude of the unequilibrated normal matrix by about
twenty-four orders of magnitude and recovers the numerical accuracy; a
further reduction to $10^{-15}$ produces little additional change. This
experiment is used only as a diagnostic of the raw system. The proposed
method applies the row equilibration in
\eqref{eq:row-equilibration} to all residual equations and introduces no
problem-dependent scaling parameter.

\paragraph{Parameter sensitivity.}
Additional tests with $\beta^+/\beta^-$ ranging from $10^{-8}$ to $10^8$
and with $M=20,\ldots,400$ random features are collected in
Figure~\ref{fig:circular_high_contrast_sensitivity} of
Appendix~\ref{app:circular_diagnostics}. Over the tested coefficient ratios,
no systematic deterioration of the relative error is observed as the
contrast magnitude increases, although exchanging the high-diffusivity phase
changes the attained error level. The feature-space study shows appreciable
dependence on the activation function and the number of features for small $M$,
while sufficiently rich sine and tanh spaces reach comparable error levels for
this benchmark. These
observations are specific to the present fixed discretization and are not
used to rank the two feature families.

\paragraph{Algebraic role of the Schur reduction.}
We next compare the proposed Schur-reduced formulation with an unreduced
full least-squares formulation (Full LS) for the same circular interface
with $(\beta^+,\beta^-)=(10^4,1)$ and $Q_1$ finite elements. In Full LS,
the interior FEM unknowns are retained in the least-squares vector and
the corresponding Galerkin equations are included as residual rows.
The comparison therefore concerns both the organization of the linear
algebra and the treatment of the interior equations as equality
constraints or residuals. Let
\[
    \widetilde{\boldsymbol F}_I
    =\boldsymbol F_I-K_{ID}\boldsymbol g_D,
    \qquad
    R_I=
    \frac{\|K_{II}\boldsymbol U_I+K_{IC}\boldsymbol U_C-
    \widetilde{\boldsymbol F}_I\|_2}{\sqrt{n_I}}.
\]
Table~\ref{tab:q1_schur_full_ls_compact} reports representative refinement
levels; the complete five-level data are given in
Table~\ref{tab:q1_schur_full_ls} of Appendix~\ref{app:circular_diagnostics}.

\begin{table}[htbp]
\centering
\caption{Representative comparison between the Schur-reduced and full
least-squares formulations for the $Q_1$ circular interface problem with
$(\beta^+,\beta^-)=(10^4,1)$.}
\label{tab:q1_schur_full_ls_compact}
\small
\setlength{\tabcolsep}{3.8pt}
\begin{tabular}{@{}c ccc ccc@{}}
\toprule
& \multicolumn{3}{c}{Schur-reduced} & \multicolumn{3}{c}{Full LS} \\
\cmidrule(lr){2-4}\cmidrule(lr){5-7}
$h$ & $E_{L^2}^{\mathrm{rel}}$ & $R_I$ & Time (s)
& $E_{L^2}^{\mathrm{rel}}$ & $R_I$ & Time (s) \\
\midrule
$1/20$ & $2.508257\times10^{-3}$ & $7.065932\times10^{-13}$ & $1.249$
& $2.508457\times10^{-3}$ & $7.236761\times10^{-6}$ & $3.468$ \\
$1/40$ & $6.747829\times10^{-4}$ & $8.585180\times10^{-13}$ & $3.675$
& $6.748271\times10^{-4}$ & $1.076821\times10^{-6}$ & $251.660$ \\
$1/80$ & $1.748283\times10^{-4}$ & $1.021440\times10^{-12}$ & $14.569$
& $1.748322\times10^{-4}$ & $7.963211\times10^{-8}$ & $10692.936$ \\
\bottomrule
\end{tabular}
\end{table}

The two formulations give essentially the same solution errors on the
reported meshes. The hard-constraint role of the Schur reduction
is visible in $R_I$, which remains near $10^{-12}$ after reconstruction
through the sparse interior response. Full LS balances the interior
equations against the other residual blocks and produces larger
interior residuals. Thus, the experiment demonstrates improved
satisfaction of the interior discrete equations.

The timing comparison illustrates a different consequence of the same
reduction. The FEM interior coefficients are obtained through sparse solves and
reconstruction. In the present Full LS implementation, these coefficients are
handled together with the coupling and random-feature coefficients by
the dense least-squares solver. The increasing cost difference in
Table~\ref{tab:q1_schur_full_ls_compact} therefore illustrates the
algebraic distinction described in Remark~\ref{rem:schur-reduction}.
These timings are specific to the present implementations and do not
constitute an asymptotic complexity estimate.

\subsection{Star-shaped interface problems}
\label{subsec:star_interface}

We consider two elliptic interface benchmarks with five-pointed
polygonal interfaces in the computational domain $\Omega=(-1,1)^2$. The present benchmarks involve nonsmooth interface geometries but piecewise smooth solutions; singular solutions associated with interface corners are not considered here.
Let $x_1=r\cos\theta$ and $x_2=r\sin\theta$, where the angular variable $\theta$ is understood modulo $2\pi$.
For both cases, the geometric parameters are $\theta_t=\frac{\pi}{5},$ $\theta_r=\frac{\pi}{7},$ $R=\frac{6}{7},$ and $i=1,\ldots,5$. 

\paragraph{Case 1.}
The first problem is taken from
\cite{hou2010numerical,ahmad2020local,wang2021new}.
The interface $\Gamma$ is represented by $\Gamma=\left\{(r,\theta):\varphi_1(r,\theta)=0\right\}$, where $\varphi_1(r,\theta)=\rho_1(\theta)-r$, and
\[
\rho_1(\theta)=
\begin{cases}
\displaystyle
\frac{R\sin\theta_t}
{2\sin\!\left(
\theta_t/2+\theta-\theta_r
-\frac{2\pi(i-1)}{5}
\right)},
&
\displaystyle
\theta_r+\frac{(2i-2)\pi}{5}
\leq\theta<
\theta_r+\frac{(2i-1)\pi}{5},
\\[3mm]
\displaystyle
\frac{R\sin\theta_t}
{2\sin\!\left(
\theta_t/2-\theta+\theta_r
+\frac{2\pi(i-1)}{5}
\right)},
&
\displaystyle
\theta_r+\frac{(2i-3)\pi}{5}
\leq\theta<
\theta_r+\frac{(2i-2)\pi}{5}.
\end{cases}
\]
The two subdomains are defined by $\Omega^+
    =
    \left\{\boldsymbol{x}\in\Omega:
    \varphi_1(\boldsymbol{x})>0\right\},$ and
    $\Omega^-
    =
    \left\{\boldsymbol{x}\in\Omega:
    \varphi_1(\boldsymbol{x})<0\right\}.$
The variable diffusion coefficient is prescribed by $\beta^+(\boldsymbol{x})=1$, $\beta^-(\boldsymbol{x})=2+\sin(x_1+x_2)$, and the exact solution is
\[
    u^+(\boldsymbol{x})
    =
    6+\sin(2\pi x_1)\sin(2\pi x_2),
    \qquad
    u^-(\boldsymbol{x})
    =
    x_1^2+x_2^2+\sin(x_1+x_2).
\]

We first compare the proposed hybrid method with the second-order IFVE
method reported in \cite{wang2021new}.
For the hybrid method, the polynomial degree is fixed at $p=1$, so that
the FEM region is discretized using the $Q_1$ finite element space.
Table~\ref{tab:star_case1_p1_comparison} reports the absolute nodal
$L^\infty$ errors for both methods. The absolute $L^2$ errors and their
observed convergence orders are additionally reported for the proposed
method to document its $Q_1$ mesh-refinement behavior.

\begin{table*}[htbp]
    \centering
    \caption{
        Comparison of the absolute nodal $L^\infty$ errors
        $E_{L^\infty}^{\mathrm{abs}}$ and the absolute $L^2$ errors
        $E_{L^2}^{\mathrm{abs}}$ for Case~1.
        The IFVE $L^\infty$ results are taken from
        \cite{wang2021new}. The $L^2$ errors are reported only for the
        proposed method, since corresponding $L^2$ data are not
        available for the IFVE method. The cross-method comparison is
        made at matched Cartesian mesh resolution and is not normalized
        by computational cost.
    }
    \label{tab:star_case1_p1_comparison}
    \small
    \setlength{\tabcolsep}{4.5pt}
    \renewcommand{\arraystretch}{1.05}
    \begin{tabular}{@{}ccccccc@{}}
        \toprule
        $N\times N$
        & \multicolumn{2}{c}{IFVE method \cite{wang2021new}}
        & \multicolumn{4}{c}{Proposed method ($Q_1$)}
        \\
        \cmidrule(lr){2-3}
        \cmidrule(lr){4-7}
        &
        $E_{L^\infty}^{\mathrm{abs}}$
        & Order
        &
        $E_{L^\infty}^{\mathrm{abs}}$
        & Order
        &
        $E_{L^2}^{\mathrm{abs}}$
        & Order
        \\
        \midrule

        $64\times64$
        & $1.4569\times10^{-2}$ & --
        & $2.2269\times10^{-3}$ & --
        & $2.1266\times10^{-3}$ & --
        \\

        $128\times128$
        & $4.1902\times10^{-3}$ & 1.7978
        & $5.4037\times10^{-4}$ & 2.0430
        & $5.7177\times10^{-4}$ & 1.8950
        \\

        $256\times256$
        & $1.0923\times10^{-3}$ & 1.9396
        & $1.2762\times10^{-4}$ & 2.0821
        & $1.5151\times10^{-4}$ & 1.9160
        \\

        $512\times512$
        & $2.8353\times10^{-4}$ & 1.9458
        & $3.0272\times10^{-5}$ & 2.0758
        & $3.9362\times10^{-5}$ & 1.9445
        \\

        $1024\times1024$
        & $7.4064\times10^{-5}$ & 1.9367
        & $6.9871\times10^{-6}$ & 2.1152
        & $1.0024\times10^{-5}$ & 1.9733
        \\

        \bottomrule
    \end{tabular}
\end{table*}

Both methods exhibit approximately second-order convergence in the reported
$E_{L^\infty}^{\mathrm{abs}}$, while the proposed $Q_1$ approximation also
shows a stable second-order decrease in $E_{L^2}^{\mathrm{abs}}$. The
comparison is made at matched Cartesian mesh resolution and concerns
accuracy only; the two methods use different approximation spaces and
algebraic solvers, so no computational-efficiency conclusion is drawn.

\begin{figure}[htbp]
    \centering
    \begin{subfigure}[t]{0.485\textwidth}
        \centering
        \includegraphics[width=\linewidth]
        {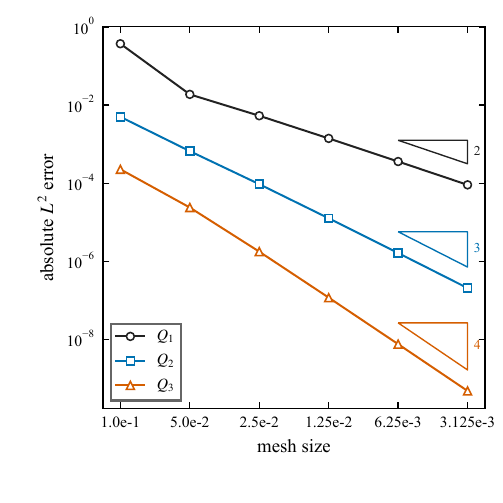}
        \caption{$E_{L^2}^{\mathrm{abs}}$.}
        \label{fig:star_case1_l2_convergence}
    \end{subfigure}
    \hfill
    \begin{subfigure}[t]{0.485\textwidth}
        \centering
        \includegraphics[width=\linewidth]
        {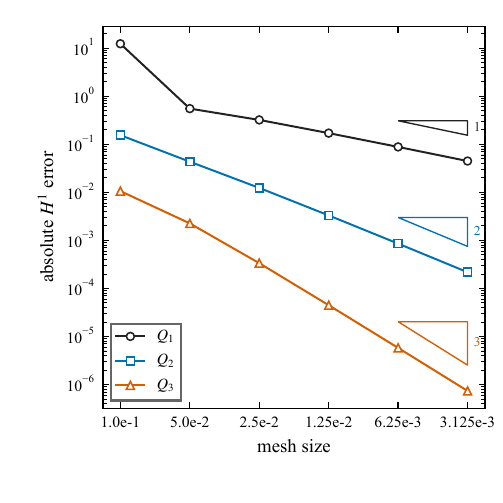}
        \caption{$E_{H^1}^{\mathrm{abs}}$.}
        \label{fig:star_case1_h1_convergence}
    \end{subfigure}
    \caption{
        Convergence of $E_{L^2}^{\mathrm{abs}}$ and
$E_{H^1}^{\mathrm{abs}}$ for Case~1 in
Section~\ref{subsec:star_interface}.
    }
    \label{fig:star_case1_convergence}
\end{figure}

We next examine high-order mesh refinement for the proposed method.
Figure~\ref{fig:star_case1_convergence} presents the absolute $L^2$
and broken $H^1$-seminorm errors for the $Q_1$--$Q_3$
discretizations.

The high-order results are consistent with the nominal $Q_1$--$Q_3$
mesh-refinement behavior. The $Q_2$ and $Q_3$ curves closely track their
reference slopes, while the $Q_1$ curves approach the expected rates after
a short pre-asymptotic regime.

\paragraph{Case 2.}
The second problem is considered in
\cite{xia2014mib,ahmad2020local,wu2022inn}.
The interface $\Gamma$ is represented by
$\Gamma=\left\{(r,\theta):\varphi_2(r,\theta)=0\right\}$, where
$\varphi_2(r,\theta)=\rho_2(\theta)-r$, and
\[
\rho_2(\theta)=
\begin{cases}
\displaystyle
\frac{R\sin(\theta_t/2)}
{\sin\!\left(
\theta_t/2+\theta-\theta_r
-\frac{2\pi(i-1)}{5}
\right)},
&
\displaystyle
\theta_r+\frac{(2i-2)\pi}{5}
\leq\theta<
\theta_r+\frac{(2i-1)\pi}{5},
\\[3mm]
\displaystyle
\frac{R\sin(\theta_t/2)}
{\sin\!\left(
\theta_t/2-\theta+\theta_r
+\frac{2\pi(i-1)}{5}
\right)},
&
\displaystyle
\theta_r+\frac{(2i-3)\pi}{5}
\leq\theta<
\theta_r+\frac{(2i-2)\pi}{5}.
\end{cases}
\]
The two subdomains are defined by
$\Omega^+
=
\left\{
\boldsymbol{x}\in\Omega:
\varphi_2(\boldsymbol{x})>0
\right\}$
and
$\Omega^-
=
\left\{
\boldsymbol{x}\in\Omega:
\varphi_2(\boldsymbol{x})<0
\right\}$.
The piecewise constant diffusion coefficient is prescribed by
$\beta^+(\boldsymbol{x})=5$ and $\beta^-(\boldsymbol{x})=3$, and the
exact solution is
\[
u^+(\boldsymbol{x})
=
7+x_1^2+x_2^2,
\qquad
u^-(\boldsymbol{x})
=
\sin(3\pi x_1)\sin(3\pi x_2).
\]

The corresponding $Q_1$--$Q_3$ convergence histories and the comparison
with the MIB--Galerkin results of \cite{xia2014mib} are reported in
Appendix~\ref{app:star_supplementary}; see
Figure~\ref{fig:star_case2_convergence} and
Table~\ref{tab:star_case2_comparison}. On the finer meshes, the $Q_1$ curves
approach second-order $L^2$ and first-order broken $H^1$ convergence, while
the $Q_2$ and $Q_3$ results follow the corresponding nominal rates. The
matched-resolution $L^\infty$ comparison is again used only as an accuracy
benchmark.

A regional decomposition of the $Q_3$ errors between the interface band and
the FEM region is also given in Appendix~\ref{app:star_supplementary},
Table~\ref{tab:star_q3_regional_relative}. In both cases, the FEM-region
errors approach fourth order in $L^2$ and third order in the broken $H^1$
seminorm, while the band errors decrease without showing a slower asymptotic
trend over the tested refinements. Because both regional domains vary with
$h$, these regional rates are used only as diagnostics on the changing
regional domains.

\subsection{A multiple-interface problem}
\label{subsec:multiple_interface}

We next consider a multiple-interface elliptic problem previously studied in
\cite{ewing1999immersed,li2025local,song2026discontinuity}.
Let $\Omega=(-1,1)^2$, $r=\sqrt{x_1^2+x_2^2}$, and
$\theta=\operatorname{atan2}(x_2,x_1)$. The three interfaces are
\[
\begin{aligned}
\Gamma_1&=\{\boldsymbol{x}\in\Omega:r=0.2\},\\
\Gamma_2&=\{\boldsymbol{x}\in\Omega:r=0.5-0.1\cos(5\theta)\},\\
\Gamma_3&=\{\boldsymbol{x}\in\Omega:r=0.8\}.
\end{aligned}
\]
They partition $\Omega$ into four nested subdomains,
\[
\begin{aligned}
\Omega_1&=\{0\le r<0.2\},&
\Omega_2&=\{0.2<r<0.5-0.1\cos(5\theta)\},\\
\Omega_3&=\{0.5-0.1\cos(5\theta)<r<0.8\},&
\Omega_4&=\{r>0.8\},
\end{aligned}
\]
with the dependence on $\boldsymbol{x}\in\Omega$ understood. The exact
solution is prescribed by
\[
u(\boldsymbol{x})
=
\begin{cases}
\cos(x_2)+1.8, & \boldsymbol{x}\in\Omega_1,\\
\exp(x_1)+1.3, & \boldsymbol{x}\in\Omega_2,\\
\sin(x_1)+0.5, & \boldsymbol{x}\in\Omega_3,\\
-x_1+\ln(x_2+2), & \boldsymbol{x}\in\Omega_4.
\end{cases}
\]
The diffusion coefficient takes the values $1,2,3,$ and $4$ in
$\Omega_1,\ldots,\Omega_4$, respectively, and the forcing, boundary data,
and interface jumps are generated consistently from the exact solution.

\begin{figure*}[t]
    \centering
    \begin{subfigure}[t]{0.485\textwidth}
        \centering
        \includegraphics[width=\linewidth]{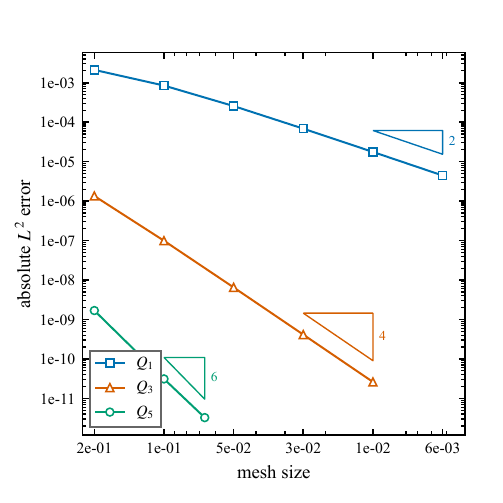}
        \caption{$E_{L^2}^{\mathrm{abs}}$.}
        \label{fig:multiple_interface_solution_convergence}
    \end{subfigure}
    \hfill
    \begin{subfigure}[t]{0.485\textwidth}
        \centering
        \includegraphics[width=\linewidth]{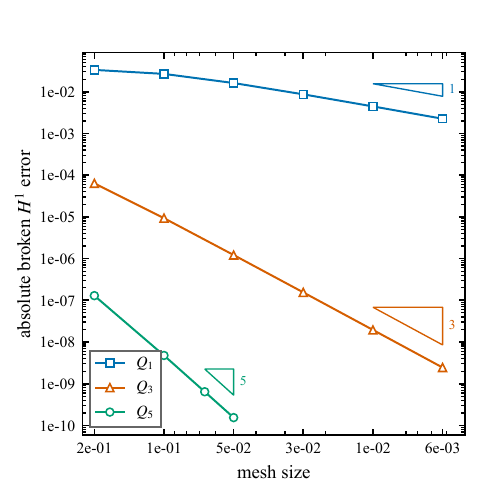}
        \caption{$E_{H^1}^{\mathrm{abs}}$.}
        \label{fig:multiple_interface_gradient_convergence}
    \end{subfigure}
    \caption{
       Convergence of $E_{L^2}^{\mathrm{abs}}$ and
$E_{H^1}^{\mathrm{abs}}$ for the multiple-interface problem in
Section~\ref{subsec:multiple_interface}.
    }
    \label{fig:multiple_interface_convergence}
\end{figure*}

This benchmark tests whether several interface bands can be treated
simultaneously while the remaining uncut cells retain the standard FEM
discretization. The same problem has also been used to assess full-domain
randomized interface solvers \cite{li2025local,song2026discontinuity}; here
the main focus is the mesh-refinement behavior of the localized hybrid
discretization. Figure~\ref{fig:multiple_interface_convergence} reports the
absolute $L^2$ and broken $H^1$ errors for the $Q_1$, $Q_3$, and $Q_5$
discretizations.

The observed rates are consistent with the nominal $Q_p$ behavior over the
tested refinement ranges: approximately orders $(2,1)$, $(4,3)$, and
$(6,5)$ in the $L^2$ and broken $H^1$ errors for $Q_1$, $Q_3$, and $Q_5$,
respectively. Thus, in this benchmark, the presence of three nested
interfaces does not visibly alter the expected $p$-dependent refinement
pattern.

Additional fixed-mesh $p$-enrichment results, regional FEM/band error
decompositions, and representative pointwise error distributions are
reported in Appendix~\ref{app:multiple_interface_supplementary}. These
supplementary diagnostics show that the band error remains below the
corresponding FEM-region error over the principal $p$-convergence range
before the computations reach the numerical accuracy floor.

\paragraph{Localized hybrid versus full-subdomain RFM approximation.}
To examine the approximation tasks and algebraic costs associated with
localization, we compare the hybrid method with full-subdomain RFM
discretizations. In these reference computations, an independent
randomized representation is assigned to each complete physical
subdomain and no FEM bulk discretization is used. This construction is
similar in spirit to the subdomain-wise randomized representations in
\cite{chi2024random,li2025local}, although the precise formulations and
coupling procedures are not identical. The same sinusoidal feature
family as in the hybrid computations is employed. For post-processing,
the full-subdomain RFM fields are also evaluated over exactly the same
union of band cells as the corresponding hybrid computation.
Table~\ref{tab:hybrid_full_rfm_comparison} reports global and band
errors, nominal total degrees of freedom, and assembly-and-solve times.
The nominal DOF counts describe the representation size and should not
be identified with the dimensions or cost of the dense least-squares
solve.

\begin{table*}[t]
\centering
\caption{
Comparison between the localized hybrid FEM--RFM formulation and
full-subdomain RFM computations for the circular interface problem in
Section~\ref{subsec:circular_high_contrast} and the multiple-interface
problem in Section~\ref{subsec:multiple_interface}.
The notation in parentheses gives the number of random-feature branches
times the number of features per branch. The band errors of all methods
are evaluated over the same hybrid band for each problem. The reported
times include assembly and solution of the discrete system.
}
\label{tab:hybrid_full_rfm_comparison}
\small
\setlength{\tabcolsep}{3.8pt}
\renewcommand{\arraystretch}{1.08}
\resizebox{\textwidth}{!}{%
\begin{tabular}{@{}llrrrrrr@{}}
\toprule
Problem
& Method
& Total DOF
& Time (s)
& $E_{L^2}^{\mathrm{abs}}$
& $E_{H^1}^{\mathrm{abs}}$
& $E_{L^2,B}^{\mathrm{abs}}$
& $E_{H^1,B}^{\mathrm{abs}}$
\\
\midrule

\multirow{3}{*}{
\shortstack[c]{
Circular interface\\
Section~\ref{subsec:circular_high_contrast}
}}
& Hybrid $Q_2$, $N=80$ ($2\times200$)
& 25,337
& 8.006
& $5.222\times10^{-7}$
& $1.355\times10^{-4}$
& $1.662\times10^{-10}$
& $1.419\times10^{-8}$
\\

& Full-subdomain RFM ($2\times200$)
& 400
& 0.113
& $4.294\times10^{-5}$
& $9.629\times10^{-4}$
& $1.410\times10^{-5}$
& $7.391\times10^{-4}$
\\

& Full-subdomain RFM ($2\times12668$)
& 25,336
& 2315.896
& $5.453\times10^{-6}$
& $2.189\times10^{-4}$
& $5.169\times10^{-7}$
& $2.633\times10^{-5}$
\\

\midrule

\multirow{4}{*}{
\shortstack[c]{
Multiple interface\\
Section~\ref{subsec:multiple_interface}
}}
& Hybrid $Q_5$, $N=30$ ($6\times200$)
& 19,271
& 25.528
& $3.238\times10^{-12}$
& $6.431\times10^{-10}$
& $9.352\times10^{-13}$
& $4.047\times10^{-12}$
\\

& Full-subdomain RFM ($4\times200$)
& 800
& 0.244
& $1.102\times10^{-9}$
& $2.006\times10^{-8}$
& $5.952\times10^{-10}$
& $7.803\times10^{-9}$
\\

& Full-subdomain RFM ($4\times1150$)
& 4,600
& 15.433
& $1.336\times10^{-12}$
& $2.567\times10^{-11}$
& $8.132\times10^{-13}$
& $1.342\times10^{-11}$
\\

& Full-subdomain RFM ($4\times4817$)
& 19,268
& 1130.697
& $1.220\times10^{-12}$
& $2.575\times10^{-11}$
& $8.107\times10^{-13}$
& $1.420\times10^{-11}$
\\

\bottomrule
\end{tabular}%
}
\end{table*}

At the fixed per-branch width $M=200$, the hybrid computations have
smaller band errors than the corresponding full-subdomain RFM
computations. For the circular problem, the number of random-feature
coefficients is the same, but the hybrid formulation also includes the
FEM field and its coupling conditions. The comparison therefore reflects differences in both spatial support
and coupling between the two formulations. For the
multiple-interface problem, only the per-branch width is matched,
since the decompositions contain different numbers of branches.
The results are consistent with a reduced local approximation task,
but the small full-subdomain RFM systems remain considerably cheaper.

The approximately matched-DOF computations illustrate the distinction
between representation size and algebraic solution cost. In the
full-subdomain RFM implementation, increasing the feature count directly
enlarges the dense least-squares unknown vector. In the hybrid method,
the interior FEM coefficients are determined by sparse
interior solves and reconstruction. Similar nominal DOF counts can
therefore correspond to different computational costs, as observed in
Table~\ref{tab:hybrid_full_rfm_comparison}. However, matching DOFs
does not establish which method is more economical at a prescribed
accuracy: a full-subdomain RFM may require substantially fewer
coefficients to reach that accuracy.

The multiple-interface benchmark provides such an example. With
$M=1150$ features per physical subdomain, the full-subdomain RFM
already attains smaller global errors than the reported $Q_5$ hybrid
computation at lower cost. In the circular benchmark, by contrast, the reduced-regularity point at the
origin is assigned to the FEM region, while the random features represent
the solution only in a smooth neighborhood of $\Gamma$. Thus, localization changes both the spatial extent of
the randomized representation and the solution behavior it must
approximate. Its computational value depends on this change in
approximation task together with the cost of the coupled algebraic
system. These comparisons are specific to the tested feature spaces and
implementations.

\subsection{An interface problem with an irregular outer boundary}
\label{subsec:irregular_outer_boundary}

We next test the irregular-boundary extension of
Remark~\ref{rem:irregular-boundary} using a benchmark adapted from
Example~6.8 of \cite{chen2025fft}. In polar coordinates
$x_1=r\cos\theta$, $x_2=r\sin\theta$,
\[
\Omega=\{(r,\theta):0\le r<1.1+0.1\sin(6\theta)\},
\qquad
\Gamma=\{(r,\theta):r=0.6+0.1\sin(6\theta)\}.
\]
The two physical subdomains are
$\Omega^-=\{r<0.6+0.1\sin(6\theta)\}$ and
$\Omega^+=\Omega\setminus\overline{\Omega^-}$, with
$(\beta^+,\kappa^+)=(2,1)$ and
$(\beta^-,\kappa^-)=(0.1,10)$. Here,
$\kappa|_{\Omega^s}=\kappa^s$ is piecewise constant. This extends the
constant-$\kappa$ core model by using the corresponding local value
in each volume operator, while retaining the same trace and
normal-flux coupling. The exact solution is
\[
u(\boldsymbol{x})
=
\begin{cases}
e^{x_1}\cos x_2, & \boldsymbol{x}\in\Omega^-,\\[1mm]
x_1^2+\sin(x_2^2), & \boldsymbol{x}\in\Omega^+,
\end{cases}
\]
and the source, Dirichlet data, and interface jumps are generated
consistently.

\begin{figure}[htbp]
    \centering
    \begin{subfigure}[t]{0.485\textwidth}
        \centering
        \includegraphics[width=\linewidth]{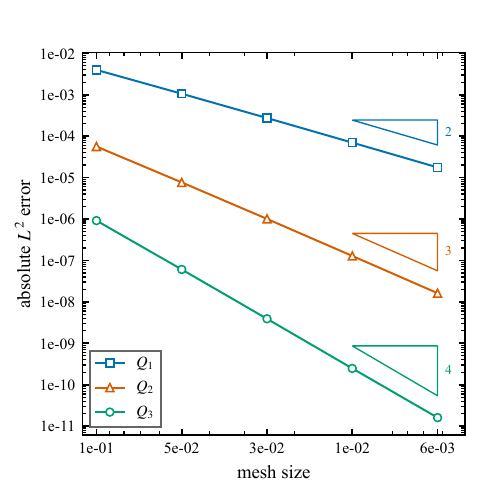}
        \caption{$E_{L^2}^{\mathrm{abs}}$.}
        \label{fig:irregular_domain_l2_convergence}
    \end{subfigure}
    \hfill
    \begin{subfigure}[t]{0.485\textwidth}
        \centering
        \includegraphics[width=\linewidth]{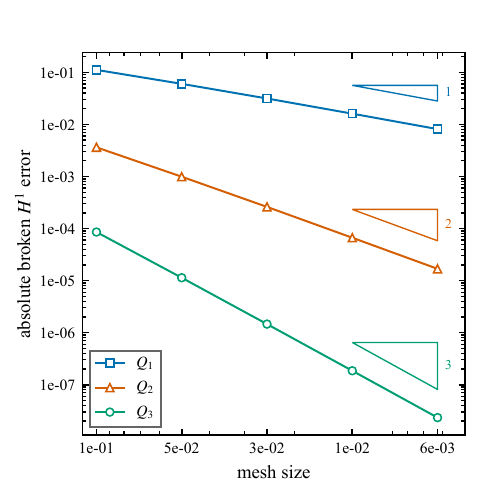}
        \caption{$E_{H^1}^{\mathrm{abs}}$.}
        \label{fig:irregular_domain_h1_convergence}
    \end{subfigure}
    \caption{
       Convergence of $E_{L^2}^{\mathrm{abs}}$ and
$E_{H^1}^{\mathrm{abs}}$ for the irregular-domain problem in
Section~\ref{subsec:irregular_outer_boundary}.
    }
    \label{fig:irregular_domain_convergence}
\end{figure}

Because neither $\Gamma$ nor $\partial\Omega$ is aligned with the Cartesian
background mesh, localized random-feature bands are used both around the
internal interface and along the irregular outer boundary, with standard
finite elements retained in the intervening uncut cells. A representative
partition is shown in Appendix~\ref{app:extension_visual_diagnostics}.
For the refinement study, $\Omega$ is embedded in $(-1.25,1.25)^2$ and
$Q_1$--$Q_3$ finite elements are used in the FEM region.
Figure~\ref{fig:irregular_domain_convergence} reports the resulting errors.

The observed rates are consistent with the nominal $Q_1$--$Q_3$ behavior:
approximately orders $(2,1)$, $(3,2)$, and $(4,3)$ in the $L^2$ and broken
$H^1$ errors, respectively. Thus, for this smooth test geometry, adding a
localized random-feature boundary band does not visibly alter the expected
mesh-refinement pattern. This is a numerical observation for the present benchmark.

\subsection{A bicontinuous composite for thermal homogenization}
\label{subsec:bicontinuous_homogenization}

We next consider a periodic thermal-homogenization problem motivated by the
bicontinuous composite studied in \cite{gehrig2025x}. Following the
Gaussian random field construction of \cite{soyarslan20183d}, we generate
an independent two-dimensional periodic realization on
$\Omega=(-1,1)^2$,
\begin{equation}
\label{eq:bicontinuous-grf}
\phi(\boldsymbol{x})
=
\sqrt{\frac{2}{2000}}
\sum_{m=1}^{2000}
\cos\!\left(
\pi\boldsymbol{k}_m\cdot\boldsymbol{x}+\theta_m
\right),
\qquad
\boldsymbol{k}_m\in\mathbb{Z}^2,\quad
|\boldsymbol{k}_m|^2=65,
\end{equation}
with independent $\theta_m\sim\mathcal U[0,2\pi)$. The zero level set
$\Gamma=\{\phi=0\}$ separates
$\Omega^-=\{\phi<0\}$ and $\Omega^+=\{\phi>0\}$, and the conductivities
are $\beta^-=0.2$ and $\beta^+=1.4$.

For the macroscopic loading $\boldsymbol{e}_x=(1,0)^T$, let $\chi$ be the
periodic corrector and define
\begin{equation}
\label{eq:complex-flux}
\boldsymbol{\xi}=\boldsymbol{e}_x+\nabla\chi,
\qquad
\boldsymbol{q}=\beta\boldsymbol{\xi}.
\end{equation}
The interface conditions are
\begin{equation}
\label{eq:complex-pde}
-\beta^\pm\Delta\chi^\pm=0
\quad\text{in }\Omega^\pm,\qquad
\jump{\chi}=0,\qquad
\jump{\beta\partial_{\boldsymbol n}\chi}
=-(\beta^+-\beta^-)n_x
\quad\text{on }\Gamma,
\end{equation}
where $n_x=\boldsymbol n\cdot\boldsymbol e_x$. The corrector and heat
flux satisfy
\begin{equation}
\label{eq:bicontinuous-periodic}
\begin{aligned}
\chi(1,x_2)&=\chi(-1,x_2),
& q_x(1,x_2)&=q_x(-1,x_2),\\
\chi(x_1,1)&=\chi(x_1,-1),
& q_y(x_1,1)&=q_y(x_1,-1).
\end{aligned}
\end{equation}
The additive constant is fixed at one FEM node.

The interface is substantially more complex than in the manufactured
benchmarks, so the one-sided random-feature regions are further subdivided
according to Remark~\ref{rem:local-rfm-subdivision}. The resulting
partition contains $269$ local RFM subdomains, each with $M=50$ fixed
features, on the same $64\times64$ background mesh for the $Q_1$ and $Q_2$
computations.

\begin{figure}[htbp]
    \centering
    \includegraphics[width=0.9\textwidth]
        {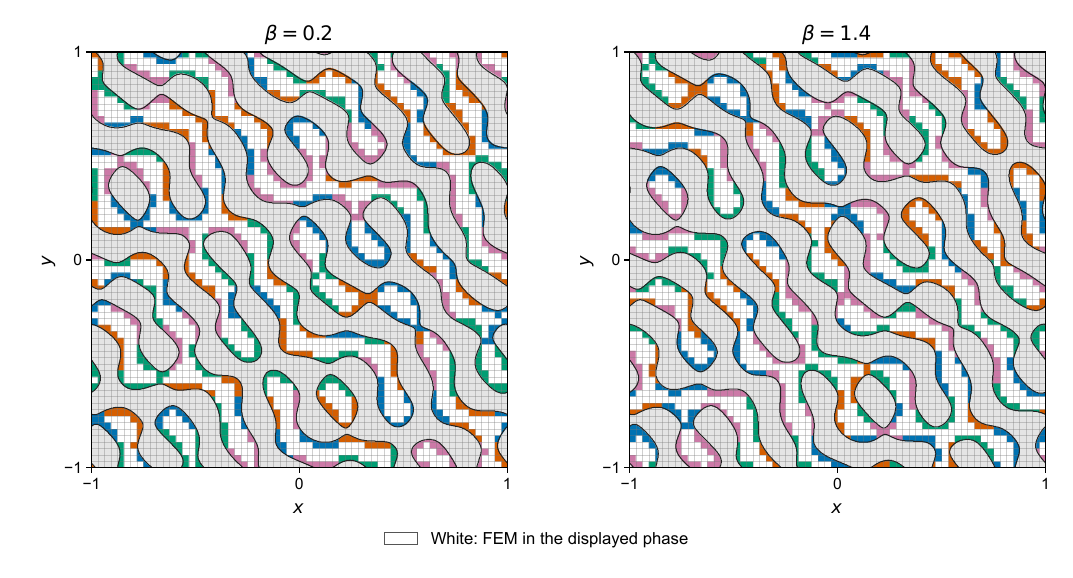}
    \caption{Hybrid partition for the bicontinuous composite.}
    \label{fig:complex-partition}
\end{figure}

\begin{figure*}[t]
    \centering
    \begin{subfigure}[t]{0.33\textwidth}
        \centering
        \includegraphics[width=\linewidth]
            {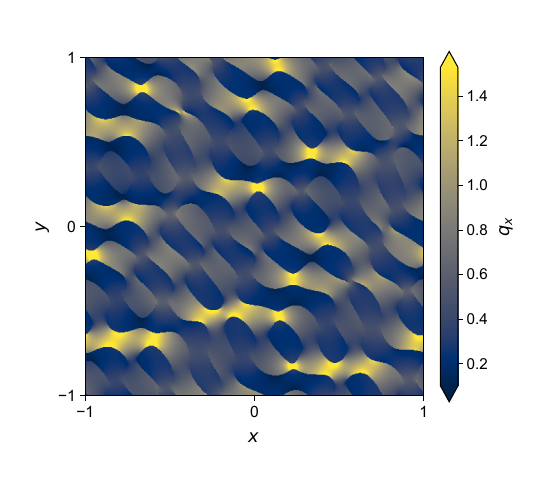}
        \caption{Reference, $N=4096$.}
    \end{subfigure}
    \hfill
    \begin{subfigure}[t]{0.33\textwidth}
        \centering
        \includegraphics[width=\linewidth]
            {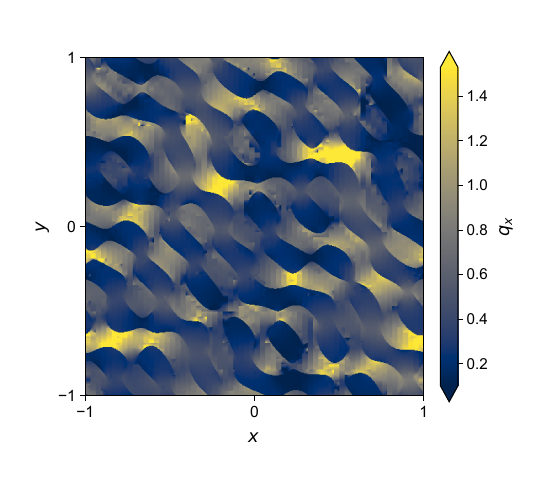}
        \caption{Hybrid $Q_1$, $N=64$.}
    \end{subfigure}
    \hfill
    \begin{subfigure}[t]{0.33\textwidth}
        \centering
        \includegraphics[width=\linewidth]
            {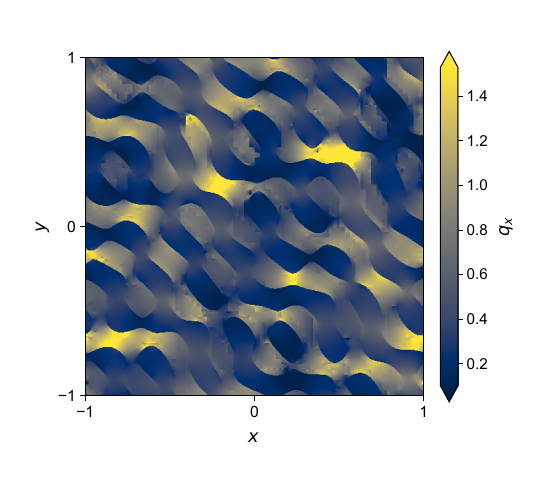}
        \caption{Hybrid $Q_2$, $N=64$.}
    \end{subfigure}
    \caption{Local heat-flux component $q_x$.}
    \label{fig:complex-qx}
\end{figure*}

Since no analytical solution is available, an independent pixel-based
$Q_1$ finite element solution on a $4096\times4096$ grid is used as the
reference. Figure~\ref{fig:complex-qx} compares the local heat-flux
component
\[
q_x=\beta(1+\partial_{x_1}\chi)
\]
for the reference and the two hybrid approximations. The corresponding
negative-temperature-gradient field $\xi_x=1+\partial_{x_1}\chi$ is shown
in Appendix~\ref{app:extension_visual_diagnostics}.

Both hybrid solutions reproduce the dominant phase-dependent structure of
the reference heat-flux field on the much coarser $64\times64$ background
mesh. For a quantitative comparison we use the effective flux
\[
\overline q_x
=
\frac{1}{|\Omega|}
\int_\Omega q_x(\boldsymbol{x})\,\mathrm d\boldsymbol{x},
\qquad
e_{\mathrm{eff}}
=
\frac{|\overline q_{x,h}-\overline q_{x,\mathrm{ref}}|}
{|\overline q_{x,\mathrm{ref}}|}.
\]

\begin{table}[htbp]
    \centering
    \caption{Effective heat flux for the bicontinuous composite.}
    \label{tab:complex-effective}
    \small
    \begin{tabular}{lccc}
        \toprule
        Method
        & $N$
        & $\overline q_x$
        & $e_{\mathrm{eff}}$ (\%)
        \\
        \midrule
        Pixel-based $Q_1$ reference
        & 4096 & 0.521364894 & --
        \\
        Hybrid $Q_1$, $M=50$
        & 64 & 0.513533233 & 1.5021
        \\
        Hybrid $Q_2$, $M=50$
        & 64 & 0.516146385 & 1.0009
        \\
        \bottomrule
    \end{tabular}
\end{table}

The relative effective-flux errors are approximately $1.50\%$ and $1.00\%$
for the $Q_1$ and $Q_2$ hybrid solutions, respectively. This experiment
therefore provides an application-oriented test of the local subdivision
strategy: randomized approximations are distributed along the complicated
physical interface, while the remaining uncut cells retain the standard
finite element discretization.

\subsection{Three-dimensional interface problems}
\label{subsec:three_dimensional_interface}

We next examine two three-dimensional tests on
$\Omega=(-1.2,1.2)^3$.

\paragraph{Case 1.}
The first interface is
$\Gamma=\{\boldsymbol{x}\in\Omega:\phi_1(\boldsymbol{x})=0\}$, where
\[
\begin{aligned}
\phi_1(\boldsymbol{x})
={}&
\left[(x_1^2+x_2^2-1)^2+x_3^2\right]
\left[(x_1^2+x_3^2-1)^2+x_2^2\right]\\
&\times
\left[(x_2^2+x_3^2-1)^2+x_1^2\right]
-(0.075)^2\left[1+3(x_1^2+x_2^2+x_3^2)\right].
\end{aligned}
\]
With $\Omega^-=\{\phi_1<0\}$ and $\Omega^+=\{\phi_1>0\}$, the diffusion
coefficient is $1$ and $100$, respectively, and the exact solution is
$u=\phi_1/\beta$. Hence $f=-\Delta\phi_1$ and both interface jumps are
homogeneous. This case combines a nonconvex three-dimensional interface
with a $100:1$ coefficient contrast. The cut-cell band is used with
$M=1500$ fixed random features on each physical side.

\paragraph{Case 2.}
The second interface is the heart-shaped surface
$\Gamma=\{\boldsymbol{x}\in\Omega:\phi_2(\boldsymbol{x})=0\}$, with
\[
\phi_2(\boldsymbol{x})
=
\left(
x_1^2+\frac94x_2^2+x_3^2-1
\right)^3
-x_1^2x_3^3
-\frac9{80}x_2^2x_3^3.
\]
Here $\beta=1$ in $\Omega^-=\{\phi_2<0\}$ and $\beta=10$ in
$\Omega^+=\{\phi_2>0\}$, while
\[
u(\boldsymbol{x})
=
\begin{cases}
\sin(r^2), & \boldsymbol{x}\in\Omega^-,\\
\cos(r^2), & \boldsymbol{x}\in\Omega^+,
\end{cases}
\qquad
r^2=x_1^2+x_2^2+x_3^2.
\]
The remaining data are generated consistently, so both the solution and
flux jumps are generally nonhomogeneous. The cut-cell band is again used,
with $M=1500$ fixed random features on each side.

\begin{figure*}[t]
    \centering
    \begin{subfigure}[t]{0.485\textwidth}
        \centering
        \includegraphics[width=\linewidth]
        {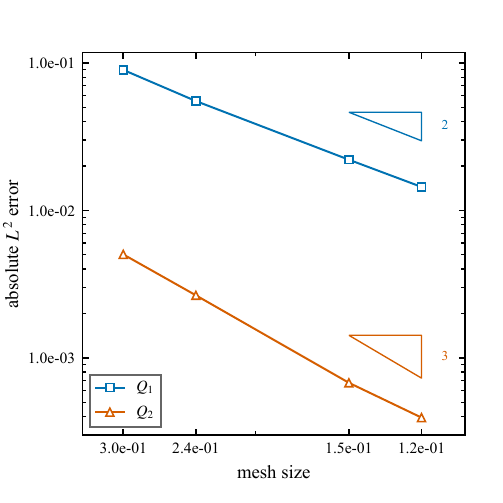}
        \caption{Case~1: absolute $L^2$ error.}
        \label{fig:three_dimensional_case1_l2}
    \end{subfigure}
    \hfill
    \begin{subfigure}[t]{0.485\textwidth}
        \centering
        \includegraphics[width=\linewidth]
        {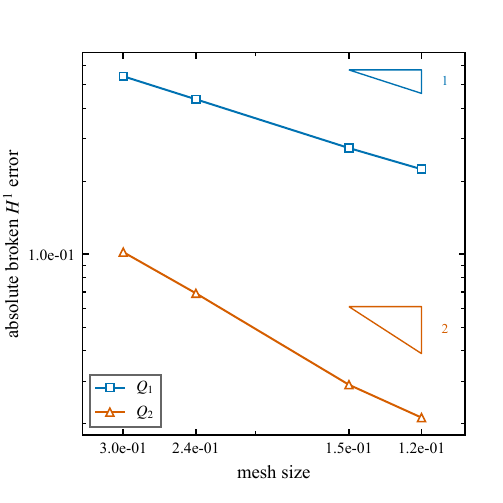}
        \caption{Case~1: absolute broken $H^1$ error.}
        \label{fig:three_dimensional_case1_h1}
    \end{subfigure}

    \begin{subfigure}[t]{0.485\textwidth}
        \centering
        \includegraphics[width=\linewidth]
        {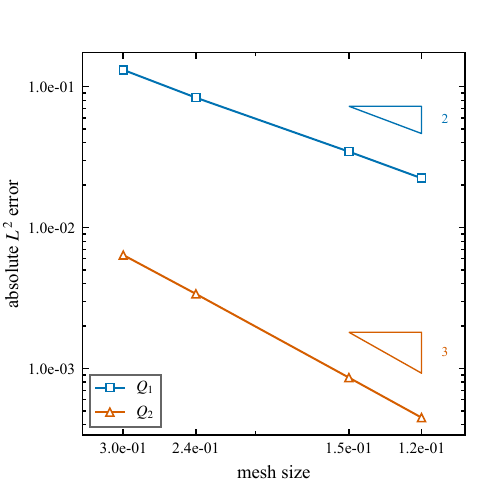}
        \caption{Case~2: absolute $L^2$ error.}
        \label{fig:three_dimensional_case2_l2}
    \end{subfigure}
    \hfill
    \begin{subfigure}[t]{0.485\textwidth}
        \centering
        \includegraphics[width=\linewidth]
        {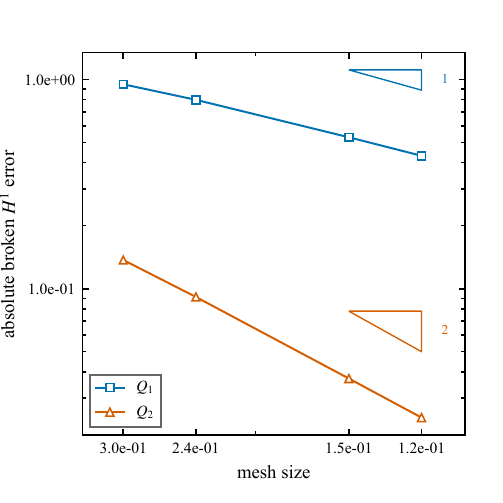}
        \caption{Case~2: absolute broken $H^1$ error.}
        \label{fig:three_dimensional_case2_h1}
    \end{subfigure}
    \caption{Mesh-refinement results for the two three-dimensional interface
    problems in Section~\ref{subsec:three_dimensional_interface}.}
    \label{fig:three_dimensional_convergence}
\end{figure*}

Figure~\ref{fig:three_dimensional_convergence} collects the $Q_1$ and $Q_2$
mesh-refinement results for both cases.

For both geometries, the $Q_1$ results are close to second order in $L^2$
and first order in the broken $H^1$ seminorm, while the $Q_2$ results track
the corresponding third- and second-order rates over the tested meshes.
Thus the expected low-order $Q_p$ refinement behavior is retained in these
two three-dimensional tests, despite the nontrivial geometry, coefficient
contrast, and nonhomogeneous jumps. Representative $Q_3$ pointwise trace
errors on the physical interfaces are provided as supplementary diagnostics
in Appendix~\ref{app:extension_visual_diagnostics}; they are not used to
infer a $Q_3$ convergence rate.

\subsection{Extensions of the proposed framework}
\label{subsec:numerical_extensions}

The preceding experiments concern the stationary linear elliptic setting.
We conclude the numerical section with two extensions of the same localized
spatial framework. These tests illustrate how the proposed framework can be
adapted to the two settings considered below.

\subsubsection{A quasilinear elliptic interface problem}
\label{subsec:quasilinear_circular_interface}

We first consider the quasilinear benchmark studied in
\cite{adjerid2025immersed,lang2026perturbation}. Let
$\Omega=(0,1)^2$ and
\[
\varphi(\boldsymbol{x})
=
(x_1-0.719)^2+(x_2-0.713)^2-0.17^2,
\qquad
\Gamma=\{\varphi=0\}.
\]
The subdomains are $\Omega^-=\{\varphi<0\}$ and
$\Omega^+=\{\varphi>0\}$. The solution-dependent diffusion coefficient is
\[
\beta(\boldsymbol{x},u)
=
\begin{cases}
\displaystyle
10\left(5+\frac{1}{\sqrt{1+100u^2}}\right),
& \boldsymbol{x}\in\Omega^-,\\[3mm]
\displaystyle
5+\frac{1}{\sqrt{1+u^2}},
& \boldsymbol{x}\in\Omega^+,
\end{cases}
\]
so $50\le\beta^-\le60$ and $5\le\beta^+\le6$. The exact solution is
\[
u(\boldsymbol{x})
=
\begin{cases}
\displaystyle
\frac1{10}\cos(x_1+3x_2)e^{x_1}\varphi(\boldsymbol{x}),
& \boldsymbol{x}\in\Omega^-,\\[2mm]
\displaystyle
\cos(x_1+3x_2)e^{x_1}\varphi(\boldsymbol{x}),
& \boldsymbol{x}\in\Omega^+,
\end{cases}
\]
with the remaining data generated consistently.

\begin{figure}[htbp]
    \centering
    \begin{subfigure}[t]{0.485\textwidth}
        \centering
        \includegraphics[width=\linewidth]
        {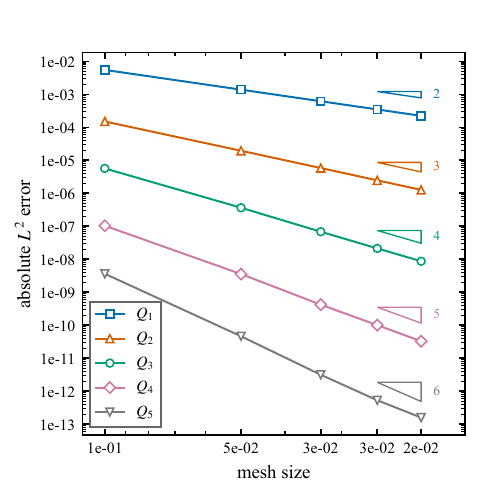}
        \caption{Absolute $L^2$ error.}
        \label{fig:quasilinear_circular_solution_convergence}
    \end{subfigure}
    \hfill
    \begin{subfigure}[t]{0.485\textwidth}
        \centering
        \includegraphics[width=\linewidth]
        {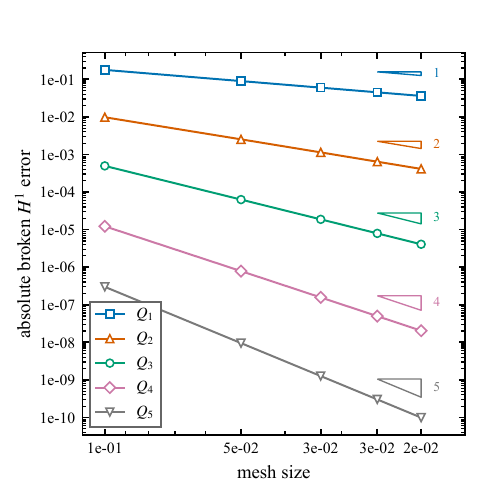}
        \caption{Absolute broken $H^1$-seminorm error.}
        \label{fig:quasilinear_circular_gradient_convergence}
    \end{subfigure}
    \caption{
        Convergence of $E_{L^2}^{\mathrm{abs}}$ and
$E_{H^1}^{\mathrm{abs}}$ for the quasilinear circular interface problem in
Section~\ref{subsec:quasilinear_circular_interface}.
    }
    \label{fig:quasilinear_circular_convergence}
\end{figure}

The $Q_1$--$Q_5$ sequences follow the corresponding nominal
$p$-dependent rates over the displayed refinement ranges, with no visible
loss of order in the high-order computations. For this benchmark, making
the diffusion coefficient solution-dependent therefore does not alter the
systematic mesh-refinement pattern observed in the linear tests. This is a
numerical observation for the present quasilinear problem.

\subsubsection{A parabolic problem with a moving interface}
\label{subsec:moving_interface_parabolic}

We next consider the moving-interface benchmark of
\cite{guo2021solving,chi2024random}. Let $\Omega=(-1,1)^2$,
$R=\pi/6$, and
\[
\rho(\boldsymbol{x},t)
=
\sqrt{
\left(x_1-0.3\cos(\pi t)\right)^2+
\left(x_2-0.3\sin(\pi t)\right)^2
}.
\]
The interface is $\Gamma(t)=\{\rho=R\}$, with
$\Omega^-(t)=\{\rho<R\}$, $\Omega^+(t)=\{\rho>R\}$,
$\beta^-=1$, and $\beta^+=10$. The exact solution is
\[
u(\boldsymbol{x},t)
=
\begin{cases}
\displaystyle
\frac{\rho(\boldsymbol{x},t)^5}{R\beta^-},
& \boldsymbol{x}\in\Omega^-(t),\\[3mm]
\displaystyle
\frac{\rho(\boldsymbol{x},t)^5}{R\beta^+}
+
R^4\left(\frac1{\beta^-}-\frac1{\beta^+}\right),
& \boldsymbol{x}\in\Omega^+(t).
\end{cases}
\]
The source, initial and boundary data, and interface jumps are generated
consistently. A representative moving-band partition is shown in
Appendix~\ref{app:extension_visual_diagnostics}.

The $Q_1$ spatial discretization is combined with BDF1 and the $Q_2$
discretization with BDF2. For spatial resolution $N$, we take
$\Delta t=N^{-2}$ and integrate to $T=1$.
Table~\ref{tab:convergence} reports the relative errors.

\begin{table}[htbp]
\centering
\caption{Relative errors and convergence orders for the moving-circle
benchmark at $T=1$. The time step is $\Delta t=N^{-2}$.}
\label{tab:convergence}
\small
\resizebox{\linewidth}{!}{%
\begin{tabular}{c cccc cccc}
\toprule
& \multicolumn{4}{c}{$Q_1$ + BDF1}
& \multicolumn{4}{c}{$Q_2$ + BDF2} \\
\cmidrule(lr){2-5} \cmidrule(lr){6-9}
$N$
& \multicolumn{2}{c}{Relative $L^2$ error}
& \multicolumn{2}{c}{Relative broken $H^1$ error}
& \multicolumn{2}{c}{Relative $L^2$ error}
& \multicolumn{2}{c}{Relative broken $H^1$ error} \\
\cmidrule(lr){2-3} \cmidrule(lr){4-5}
\cmidrule(lr){6-7} \cmidrule(lr){8-9}
& error & order
& error & order
& error & order
& error & order \\
\midrule
10
& $5.099\times10^{-2}$ & --
& $1.603\times10^{-1}$ & --
& $1.169\times10^{-3}$ & --
& $1.069\times10^{-2}$ & -- \\
20
& $1.320\times10^{-2}$ & 1.95
& $8.182\times10^{-2}$ & 0.97
& $1.615\times10^{-4}$ & 2.86
& $2.954\times10^{-3}$ & 1.86 \\
30
& $5.943\times10^{-3}$ & 1.97
& $5.493\times10^{-2}$ & 0.98
& $4.940\times10^{-5}$ & 2.92
& $1.355\times10^{-3}$ & 1.92 \\
40
& $3.371\times10^{-3}$ & 1.97
& $4.141\times10^{-2}$ & 0.98
& $2.129\times10^{-5}$ & 2.93
& $7.786\times10^{-4}$ & 1.93 \\
50
& $2.167\times10^{-3}$ & 1.98
& $3.317\times10^{-2}$ & 0.99
& $1.095\times10^{-5}$ & 2.98
& $5.006\times10^{-4}$ & 1.98 \\
\bottomrule
\end{tabular}%
}
\end{table}

The $Q_1$--BDF1 results are approximately second order in relative $L^2$
and first order in relative broken $H^1$, while the $Q_2$--BDF2 results are
approximately third and second order, respectively. With
$\Delta t=N^{-2}$, the temporal discretization does not reduce the nominal
spatial rates over the reported range. The experiment therefore shows that
the localized spatial construction can be combined with standard implicit
time stepping for this prescribed moving-interface problem.

\section{Conclusion}

We developed a localized hybrid FEM--RFM
for elliptic interface problems on unfitted background meshes.
The formulation assigns standard $Q_p$ finite elements to uncut cells
and independent random-feature representations to the interface band.
The physical jump conditions on $\Gamma$ are separated from the
trace and variational flux coupling on $\Gamma_c$. This decomposition is complemented by a Schur reduction that retains the
interior Galerkin equations as equality constraints and restricts the
dense least-squares unknown vector to the coupling nodal values and
random-feature coefficients. Thus, spatial localization of the randomized
approximation is accompanied by a reduction of the dense coupled
system that preserves the standard FEM interior discretization.

For the tested discretizations and refinement ranges, the observed
$L^2$ and broken $H^1$ errors exhibit the nominal $p$-dependent
behavior on benchmarks with sufficiently regular subdomain solutions.
The algebraic diagnostics show the importance of residual
equilibration in the high-contrast tests and confirm the preservation
of the interior FEM equations after Schur reconstruction. The Schur
comparison also shows that nearly identical solution errors can be
obtained with different costs and interior residual levels. The
full-subdomain RFM comparisons distinguish a change in approximation
task from a change in algebraic cost: localization can reduce the observed band approximation error at a fixed per-branch feature count in the tested configurations, but a smooth
full-subdomain solution may be represented more economically by a
smaller purely randomized discretization. No uniform efficiency
advantage is inferred from matched nominal DOF counts.

Additional tests with an irregular outer boundary, a quasilinear elliptic equation, and a prescribed moving interface illustrate the reuse of the trace and variational flux coupling under changes in geometry and governing equation. The periodic bicontinuous composite further provides an application-oriented thermal-homogenization test of the localized band representation and its local subdivision strategy. Theoretical stability and error analyses of the coupled discretization
and extensions to fully coupled free-boundary problems remain natural
directions for future work.

\section*{Acknowledgments}
This work was funded by the National Natural Science Foundation of China (grant No. 12571431).
%Bibliography
\bibliographystyle{unsrt}  
\bibliography{references}  
\newpage
\appendix

\section{Additional numerical results}

\subsection{Circular-interface diagnostics}
\label{app:circular_diagnostics}

This subsection collects additional diagnostics for the circular-interface
benchmark in Section~\ref{subsec:circular_high_contrast}, including residual
scaling, parameter sensitivity, the effect of Schur reduction, and a comparison
of strong- and weak-form residuals in the random-feature band.

\begin{table}[htbp]
    \centering
    \caption{
        Effect of scaling the physical flux-jump residual in the
        unequilibrated reduced least-squares system for the $Q_2$
        discretization with $N=128$ and $h=1/64$.
        Row equilibration is deliberately disabled in this diagnostic
        experiment. The quantity
        $\lambda_{\max}(A_s^T A_s)$ measures the spectral magnitude
        of the associated normal matrix.
    }
    \label{tab:circular_high_contrast_scaling}
    \small
    \setlength{\tabcolsep}{5.5pt}
    \renewcommand{\arraystretch}{1.08}
    \begin{tabular}{@{}ccccc@{}}
        \toprule
        $(\beta^+,\beta^-)$
        & $s$
        & $E_{L^2}^{\mathrm{abs}}$
        & $E_{H^1}^{\mathrm{abs}}$
        & $\lambda_{\max}(A_s^{T}A_s)$
        \\
        \midrule

        $(10^{12},1)$
        & $1$
        & $1.4282\times10^{-1}$
        & $3.6015\times10^{-1}$
        & $1.8220\times10^{33}$
        \\

        $(10^{12},1)$
        & $10^{-6}$
        & $6.0529\times10^{-6}$
        & $8.4556\times10^{-5}$
        & $1.8220\times10^{21}$
        \\

        $(10^{12},1)$
        & $10^{-12}$
        & $1.2920\times10^{-7}$
        & $5.3594\times10^{-5}$
        & $1.8221\times10^{9}$
        \\

        $(10^{12},1)$
        & $10^{-15}$
        & $1.3026\times10^{-7}$
        & $5.3594\times10^{-5}$
        & $1.5144\times10^{9}$
        \\

        \midrule

        $(1,10^{12})$
        & $1$
        & $8.5346\times10^{-3}$
        & $1.4864\times10^{-1}$
        & $9.0438\times10^{31}$
        \\

        $(1,10^{12})$
        & $10^{-6}$
        & $2.8371\times10^{-6}$
        & $1.1259\times10^{-4}$
        & $9.0438\times10^{19}$
        \\

        $(1,10^{12})$
        & $10^{-12}$
        & $2.9770\times10^{-7}$
        & $1.1239\times10^{-4}$
        & $9.0697\times10^{7}$
        \\

        $(1,10^{12})$
        & $10^{-15}$
        & $2.9770\times10^{-7}$
        & $1.1239\times10^{-4}$
        & $8.8694\times10^{7}$
        \\

        \bottomrule
    \end{tabular}
\end{table}

\begin{figure}[htbp]
    \centering
    \begin{subfigure}[t]{0.48\textwidth}
        \centering
        \includegraphics[width=\linewidth]{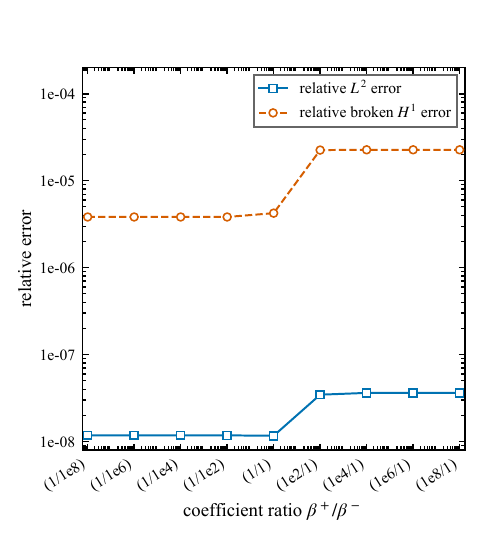}
        \caption{Relative errors versus $\beta^+/\beta^-$.}
        \label{fig:circular_coefficient_ratio_sensitivity}
    \end{subfigure}
    \hfill
    \begin{subfigure}[t]{0.48\textwidth}
        \centering
        \includegraphics[width=\linewidth]{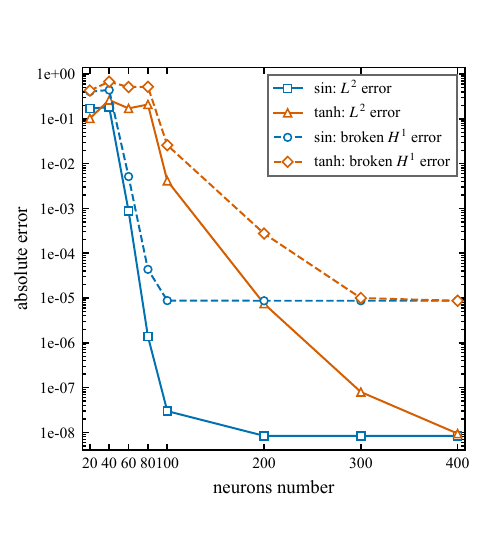}
        \caption{Absolute errors versus $M$ for sine and tanh features.}
        \label{fig:circular_feature_sensitivity}
    \end{subfigure}
    \caption{Parameter sensitivity for the circular interface problem.}
    \label{fig:circular_high_contrast_sensitivity}
\end{figure}

\begin{table}[htbp]
    \centering
    \caption{
        Comparison between the Schur-reduced and full least-squares
        formulations for the $Q_1$ circular interface problem with
        $(\beta^+,\beta^-)=(10^4,1)$. Here, $R_I$ denotes the
        normalized residual of the interior FEM equations. The reported
        times refer to the present implementation and computational
        platform.
    }
    \label{tab:q1_schur_full_ls}
    \small
    \setlength{\tabcolsep}{4.2pt}
    \renewcommand{\arraystretch}{1.05}
    \begin{tabular}{@{}ccccccccc@{}}
        \toprule
        & \multicolumn{4}{c}{Schur-reduced}
        & \multicolumn{4}{c}{Full LS}
        \\
        \cmidrule(lr){2-5}
        \cmidrule(lr){6-9}
        $h$
        & $E_{L^2}^{\mathrm{rel}}$
        & Order
        & $R_I$
        & Time (s)
        & $E_{L^2}^{\mathrm{rel}}$
        & Order
        & $R_I$
        & Time (s)
        \\
        \midrule

        $1/5$
        & $2.884753\times10^{-2}$ & --
        & $4.528737\times10^{-13}$ & $0.389$
        & $2.884934\times10^{-2}$ & --
        & $5.163037\times10^{-4}$ & $0.396$
        \\

        $1/10$
        & $8.312536\times10^{-3}$ & $1.795$
        & $5.956664\times10^{-13}$ & $0.510$
        & $8.313331\times10^{-3}$ & $1.795$
        & $5.951111\times10^{-4}$ & $0.804$
        \\

        $1/20$
        & $2.508257\times10^{-3}$ & $1.729$
        & $7.065932\times10^{-13}$ & $1.249$
        & $2.508457\times10^{-3}$ & $1.729$
        & $7.236761\times10^{-6}$ & $3.468$
        \\

        $1/40$
        & $6.747829\times10^{-4}$ & $1.894$
        & $8.585180\times10^{-13}$ & $3.675$
        & $6.748271\times10^{-4}$ & $1.894$
        & $1.076821\times10^{-6}$ & $251.660$
        \\

        $1/80$
        & $1.748283\times10^{-4}$ & $1.948$
        & $1.021440\times10^{-12}$ & $14.569$
        & $1.748322\times10^{-4}$ & $1.949$
        & $7.963211\times10^{-8}$ & $10692.936$
        \\

        \bottomrule
    \end{tabular}
\end{table}

\label{app:weak_elm_fem}

The proposed hybrid formulation employs strong-form residuals for the
random-feature approximation in the interface band. To examine the
influence of this choice, we supplement the circular interface benchmark
in Section~\ref{subsec:circular_high_contrast} with a weak-form
random-feature variant.

Following the random-feature Galerkin formulation in
\cite{sun2024local}, the band PDE residuals are replaced by weak residuals
obtained through integration by parts, where the same fixed random
features are used as both trial and test functions. The corresponding
boundary flux contributions are retained consistently. All other
components, including the one-sided random-feature representations,
physical interface jump conditions, FEM--RFM coupling on $\Gamma_c$,
row equilibration, and Schur reduction, remain unchanged.

This experiment is designed to assess the sensitivity of the proposed
hybrid discretization to the choice of residual formulation in the
localized random-feature region, rather than to compare two independent
numerical methods. The same random realization and the same number of
sine features ($M=200$ per physical side) are used for both formulations.
The coefficient contrast is fixed as
$(\beta^+,\beta^-)=(10^4,1)$.
The solution errors are measured in the global $L^2$ norm and the broken
$H^1$ seminorm defined in Section~\ref{sec:numerical_examples}.

\begin{table}[htbp]
\centering
\caption{Errors obtained with strong- and weak-form band residuals for
the $Q_1$ hybrid FEM--RFM discretization of the circular interface problem
with $(\beta^+,\beta^-)=(10^4,1)$.}
\label{tab:weak_circle_q1}
\small
\setlength{\tabcolsep}{5pt}
\begin{tabular}{@{}cc cc cc@{}}
\toprule
& & \multicolumn{2}{c}{Strong-form residual}
  & \multicolumn{2}{c}{Weak-form residual} \\
\cmidrule(lr){3-4}\cmidrule(lr){5-6}
$N$ & $h$
& $E_{L^2}^{\mathrm{abs}}$ & $E_{H^1}^{\mathrm{abs}}$
& $E_{L^2}^{\mathrm{abs}}$ & $E_{H^1}^{\mathrm{abs}}$ \\
\midrule
10  & $1/5$
& $6.662931\times10^{-3}$ & $7.266344\times10^{-2}$
& $6.691057\times10^{-3}$ & $7.267362\times10^{-2}$ \\
20  & $1/10$
& $1.919954\times10^{-3}$ & $4.229310\times10^{-2}$
& $1.921696\times10^{-3}$ & $4.230333\times10^{-2}$ \\
40  & $1/20$
& $5.793360\times10^{-4}$ & $2.560420\times10^{-2}$
& $5.798916\times10^{-4}$ & $2.560631\times10^{-2}$ \\
80  & $1/40$
& $1.558556\times10^{-4}$ & $1.378860\times10^{-2}$
& $1.559356\times10^{-4}$ & $1.378862\times10^{-2}$ \\
160 & $1/80$
& $4.038036\times10^{-5}$ & $7.147168\times10^{-3}$
& $4.038732\times10^{-5}$ & $7.147169\times10^{-3}$ \\
\bottomrule
\end{tabular}
\end{table}

\begin{table}[htbp]
\centering
\caption{Errors obtained with strong- and weak-form band residuals for
the $Q_2$ hybrid FEM--RFM discretization of the circular interface problem
with $(\beta^+,\beta^-)=(10^4,1)$.}
\label{tab:weak_circle_q2}
\small
\setlength{\tabcolsep}{5pt}
\begin{tabular}{@{}cc cc cc@{}}
\toprule
& & \multicolumn{2}{c}{Strong-form residual}
  & \multicolumn{2}{c}{Weak-form residual} \\
\cmidrule(lr){3-4}\cmidrule(lr){5-6}
$N$ & $h$
& $E_{L^2}^{\mathrm{abs}}$ & $E_{H^1}^{\mathrm{abs}}$
& $E_{L^2}^{\mathrm{abs}}$ & $E_{H^1}^{\mathrm{abs}}$ \\
\midrule
10  & $1/5$
& $2.119067\times10^{-4}$ & $6.973253\times10^{-3}$
& $2.121400\times10^{-4}$ & $6.988122\times10^{-3}$ \\
20  & $1/10$
& $2.917236\times10^{-5}$ & $1.902775\times10^{-3}$
& $3.010280\times10^{-5}$ & $1.937071\times10^{-3}$ \\
40  & $1/20$
& $4.023385\times10^{-6}$ & $5.226152\times10^{-4}$
& $4.023702\times10^{-6}$ & $5.226207\times10^{-4}$ \\
80  & $1/40$
& $5.221618\times10^{-7}$ & $1.354408\times10^{-4}$
& $5.221644\times10^{-7}$ & $1.354408\times10^{-4}$ \\
160 & $1/80$
& $6.648483\times10^{-8}$ & $3.446780\times10^{-5}$
& $6.648488\times10^{-8}$ & $3.446780\times10^{-5}$ \\
\bottomrule
\end{tabular}
\end{table}

For both $Q_1$ and $Q_2$ discretizations, the strong- and weak-form
variants produce nearly identical error levels and convergence trends
over the tested refinement sequences. The weak-form residual does not
provide a noticeable accuracy improvement for this benchmark, while the
strong-form formulation avoids the additional construction of
feature-based test matrices and boundary integral terms.

These results suggest that, for the localized treatment of the interface band
considered in this work, the strong-form residual provides an effective
realization of the random-feature approximation with a simpler residual
assembly in the present implementation. The comparison is restricted to the present benchmark and
does not imply a general superiority of one residual formulation over the
other.

\subsection{Star-shaped interface results}
\label{app:star_supplementary}

This subsection supplements the star-shaped interface experiments in
Section~\ref{subsec:star_interface}. We report the complete Case~2 refinement
and benchmark comparison, followed by the regional $Q_3$ error diagnostics
for both cases and a $Q_3$ comparison of soft and hard trace coupling for
Case~1.

\begin{figure}[htbp]
\centering
\begin{subfigure}[t]{0.485\textwidth}
\centering
\includegraphics[width=\linewidth]
{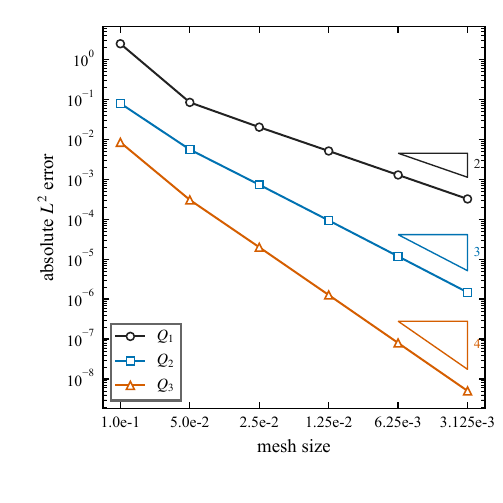}
\caption{Absolute $L^2$ error $E_{L^2}^{\mathrm{abs}}$.}
\label{fig:star_case2_l2_convergence}
\end{subfigure}
\hfill
\begin{subfigure}[t]{0.485\textwidth}
\centering
\includegraphics[width=\linewidth]
{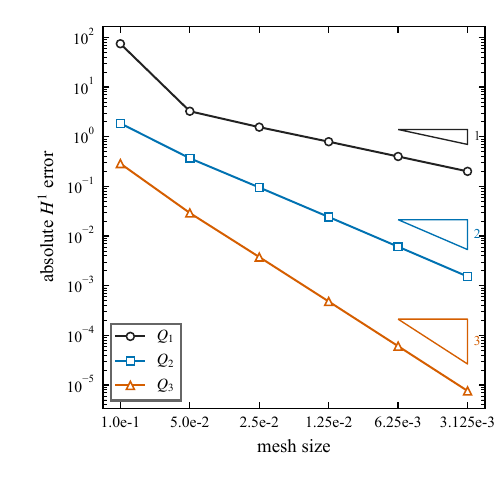}
\caption{Absolute broken $H^1$-seminorm error
$E_{H^1}^{\mathrm{abs}}$.}
\label{fig:star_case2_h1_convergence}
\end{subfigure}
\caption{
Convergence of $E_{L^2}^{\mathrm{abs}}$ and
$E_{H^1}^{\mathrm{abs}}$ for Case~2 in
Section~\ref{subsec:star_interface}.
}
\label{fig:star_case2_convergence}
\end{figure}

\begin{table}[htbp]
    \centering
    \caption{
        Comparison of the absolute $L^\infty$ errors
        $E_{L^\infty}^{\mathrm{abs}}$ for Case~2.
        The MIB--Galerkin results are taken from \cite{xia2014mib}.
        The values in parentheses denote the observed convergence
        orders. The comparison is made at matched Cartesian mesh
        resolution and is not normalized by computational cost.
    }
    \label{tab:star_case2_comparison}
    \small
    \setlength{\tabcolsep}{9pt}
    \begin{tabular}{ccc}
        \toprule
        $N\times N$
        & Proposed method ($Q_1$)
        & MIB--Galerkin \cite{xia2014mib}
        \\
        \midrule

        $20\times20$
        & $9.9959\times10^{-2}$ \; (--)
        & $6.374\times10^{-1}$ \; (--)
        \\

        $40\times40$
        & $2.5161\times10^{-2}$ \; (1.918)
        & $9.271\times10^{-2}$ \; (2.78)
        \\

        $80\times80$
        & $6.1130\times10^{-3}$ \; (2.004)
        & $2.001\times10^{-2}$ \; (2.21)
        \\

        $160\times160$
        & $1.5030\times10^{-3}$ \; (2.006)
        & $4.641\times10^{-3}$ \; (2.11)
        \\

        \bottomrule
    \end{tabular}
\end{table}

For the regional diagnostic, let $R\in\{B,F\}$ denote the interface band and
FEM regions, respectively, and define
\[
E_{L^2,R}^{\mathrm{rel}}
=\frac{\|u-u_h\|_{L^2(\Omega_h^R)}}{\|u\|_{L^2(\Omega_h^R)}},
\qquad
E_{H^1,R}^{\mathrm{rel}}
=\frac{|u-u_h|_{H^1_{\mathrm{br}}(\Omega_h^R)}}
{|u|_{H^1_{\mathrm{br}}(\Omega_h^R)}}.
\]
The broken seminorm is evaluated separately over the intersections of
$\Omega_h^R$ with the physical subdomains.

\begin{table*}[t]
\centering
\caption{
Relative regional errors for the $Q_3$ discretization of the two
five-point-star interface problems.
}
\label{tab:star_q3_regional_relative}
\small
\setlength{\tabcolsep}{3.8pt}
\renewcommand{\arraystretch}{1.05}
\resizebox{\textwidth}{!}{%
\begin{tabular}{c cc cc cc cc}
\toprule
$h$
& $E_{L^2,B}^{\mathrm{rel}}$ & Order
& $E_{L^2,F}^{\mathrm{rel}}$ & Order
& $E_{H^1,B}^{\mathrm{rel}}$ & Order
& $E_{H^1,F}^{\mathrm{rel}}$ & Order
\\
\midrule
\multicolumn{9}{c}{Case~1} \\
\midrule
0.20000
& $1.839\times10^{-5}$ & --
& $6.553\times10^{-5}$ & --
& $8.819\times10^{-4}$ & --
& $2.554\times10^{-3}$ & --
\\
0.10000
& $6.702\times10^{-7}$ & 4.78
& $5.542\times10^{-6}$ & 3.56
& $2.342\times10^{-5}$ & 5.23
& $4.890\times10^{-4}$ & 2.39
\\
0.05000
& $5.237\times10^{-8}$ & 3.68
& $3.630\times10^{-7}$ & 3.93
& $1.531\times10^{-6}$ & 3.94
& $6.784\times10^{-5}$ & 2.85
\\
0.02500
& $3.536\times10^{-9}$ & 3.89
& $2.283\times10^{-8}$ & 3.99
& $3.351\times10^{-8}$ & 5.51
& $8.721\times10^{-6}$ & 2.96
\\
0.01250
& $2.200\times10^{-10}$ & 4.01
& $1.431\times10^{-9}$ & 4.00
& $2.159\times10^{-9}$ & 3.96
& $1.105\times10^{-6}$ & 2.98
\\
0.00625
& $1.227\times10^{-11}$ & 4.16
& $8.979\times10^{-11}$ & 3.99
& $1.260\times10^{-10}$ & 4.10
& $1.392\times10^{-7}$ & 2.99
\\
\midrule
\multicolumn{9}{c}{Case~2} \\
\midrule
0.20000
& $1.202\times10^{-3}$ & --
& $1.510\times10^{-3}$ & --
& $3.266\times10^{-2}$ & --
& $2.103\times10^{-2}$ & --
\\
0.10000
& $5.112\times10^{-6}$ & 7.88
& $6.167\times10^{-5}$ & 4.61
& $8.342\times10^{-5}$ & 8.61
& $2.612\times10^{-3}$ & 3.01
\\
0.05000
& $3.203\times10^{-7}$ & 4.00
& $3.452\times10^{-6}$ & 4.16
& $3.416\times10^{-6}$ & 4.61
& $3.280\times10^{-4}$ & 2.99
\\
0.02500
& $1.902\times10^{-8}$ & 4.07
& $2.061\times10^{-7}$ & 4.07
& $2.196\times10^{-7}$ & 3.96
& $4.103\times10^{-5}$ & 3.00
\\
0.01250
& $1.105\times10^{-9}$ & 4.11
& $1.258\times10^{-8}$ & 4.03
& $1.422\times10^{-8}$ & 3.95
& $5.124\times10^{-6}$ & 3.00
\\
0.00625
& $7.742\times10^{-11}$ & 3.84
& $7.801\times10^{-10}$ & 4.01
& $1.282\times10^{-9}$ & 3.47
& $6.406\times10^{-7}$ & 3.00
\\
\bottomrule
\end{tabular}%
}
\end{table*}

The FEM-region errors approach the nominal $Q_3$ behavior in both cases.
The band errors also decrease systematically and do not exhibit a slower
asymptotic trend over the tested refinements. Since $\Omega_h^B$ and
$\Omega_h^F$ change with the background mesh, the regional rates are used only as diagnostics on the changing regional
domains.
\paragraph{Soft and hard trace coupling for Case~1.}

We further examine the treatment of trace continuity on the artificial
coupling interface $\Gamma_c$ using the $Q_3$ discretization of
Case~1. The proposed formulation uses the soft nodal trace residual
\eqref{eq:trace-residual}. For comparison, we consider two hard
alternatives that replace this residual by equality constraints while
leaving the Schur reduction, variational flux coupling, band PDE
residuals, and physical interface jump residuals unchanged.

The first alternative imposes nodal equality,
\begin{equation}
\label{eq:hard-nodal-trace}
    \boldsymbol U_C=\Phi_C\boldsymbol c.
\end{equation}
The second imposes moment equality against every free FEM trace basis
function:
\begin{equation}
\label{eq:hard-moment-trace}
    \sum_{s\in\{-,+\}}
    \int_{\Gamma_c^s}
    \left(
        u_F^h-u_B^{h,s}
    \right)
    N_{c_i}\,\mathrm ds
    =0,
    \qquad i=1,\ldots,n_C.
\end{equation}
The integrals in \eqref{eq:hard-moment-trace} are evaluated using the
same coupling quadrature as in the variational flux residual.

To isolate the effect of the trace treatment, all three formulations
use the same $Q_3$ discretization, random features, collocation sets,
row equilibration, and remaining residual equations. The computations
use $N=20,40,80,160$, $M=300$ random features and
$N_B^s=600$ band-interior collocation points on each physical side,
together with the minimal interface band defined by
\eqref{eq:band-cell-set}. No ridge regularization is used. For the hard
variants, the constrained variables are parameterized so that the
minimum-norm selection is consistent with the Euclidean norm of the
original coefficient vector $(\boldsymbol U_C,\boldsymbol c)$.

\begin{table*}[t]
\centering
\caption{
Comparison of soft and hard trace coupling for the $Q_3$
discretization of star-shaped Case~1. The errors are the absolute
$L^2$ error and absolute broken $H^1$-seminorm error.
}
\label{tab:star_case1_trace_constraints}
\small
\setlength{\tabcolsep}{5pt}
\renewcommand{\arraystretch}{1.05}
\begin{tabular}{@{}llrrrr@{}}
\toprule
Method
& $N$
& $E_{L^2}^{\mathrm{abs}}$
& Order
& $E_{H^1}^{\mathrm{abs}}$
& Order
\\
\midrule

\multirow{4}{*}{Soft}
& 20
& $2.4123\times10^{-5}$ & --
& $2.2623\times10^{-3}$ & --
\\
& 40
& $1.7957\times10^{-6}$ & 3.748
& $3.3617\times10^{-4}$ & 2.751
\\
& 80
& $1.1989\times10^{-7}$ & 3.905
& $4.4869\times10^{-5}$ & 2.905
\\
& 160
& $7.7505\times10^{-9}$ & 3.951
& $5.8013\times10^{-6}$ & 2.951
\\

\midrule

\multirow{4}{*}{Hard nodal}
& 20
& $2.8754\times10^{-5}$ & --
& $2.3040\times10^{-3}$ & --
\\
& 40
& $2.4734\times10^{-6}$ & 3.539
& $3.3686\times10^{-4}$ & 2.774
\\
& 80
& $2.3311\times10^{-7}$ & 3.407
& $4.4903\times10^{-5}$ & 2.907
\\
& 160
& $2.6490\times10^{-8}$ & 3.137
& $5.8047\times10^{-6}$ & 2.952
\\

\midrule

\multirow{4}{*}{Hard moment}
& 20
& $6.2602\times10^{-5}$ & --
& $2.9711\times10^{-3}$ & --
\\
& 40
& $6.6947\times10^{-6}$ & 3.225
& $3.4480\times10^{-4}$ & 3.107
\\
& 80
& $7.8590\times10^{-7}$ & 3.091
& $4.5362\times10^{-5}$ & 2.926
\\
& 160
& $9.8175\times10^{-8}$ & 3.001
& $5.8501\times10^{-6}$ & 2.955
\\

\bottomrule
\end{tabular}
\end{table*}

The soft formulation retains nearly fourth-order $L^2$ convergence over
the finer meshes, with an observed order of $3.95$ from $N=80$ to
$N=160$. In contrast, the corresponding fine-grid orders are $3.14$
for the hard nodal constraint and $3.00$ for the hard moment
constraint. At $N=160$, their $L^2$ errors are approximately $3.42$
and $12.67$ times the soft-coupling error, respectively. The broken
$H^1$ errors show a different behavior: all three formulations approach
an order close to three on the finest refinement.

The hard constraints are numerically satisfied to high accuracy. The
nodal mismatch is zero to the accuracy of the constrained solve, while
the relative moment residual for the hard moment formulation decreases
to $8.6\times10^{-13}$ at $N=160$; the corresponding maximum algebraic
backward residual is approximately $7.1\times10^{-17}$. The observed
difference is therefore not attributable to an inaccurate enforcement
of the hard constraints.

These results show that exact discrete trace constraints do not improve
the approximation in this test. Instead, restricting the two trace
representations by either hard constraint degrades the observed
$L^2$ convergence, whereas the broken $H^1$ convergence remains
essentially unchanged. This suggests that retaining
$\boldsymbol U_C$ and $\boldsymbol c$ as independent variables and
balancing their mismatch with the other coupled residuals can be
beneficial for the present heterogeneous FEM--RFM approximation spaces.
This interpretation is limited to the present numerical experiment and
does not constitute a general stability or convergence result.

\subsection{Multiple-interface results}
\label{app:multiple_interface_supplementary}

This subsection collects the additional fixed-mesh $p$-enrichment, regional-error,
and spatial-error results for the multiple-interface problem in
Section~\ref{subsec:multiple_interface}.

\begin{figure*}[t]
    \centering
    \begin{subfigure}[t]{0.45\textwidth}
        \centering
        \includegraphics[width=\linewidth]{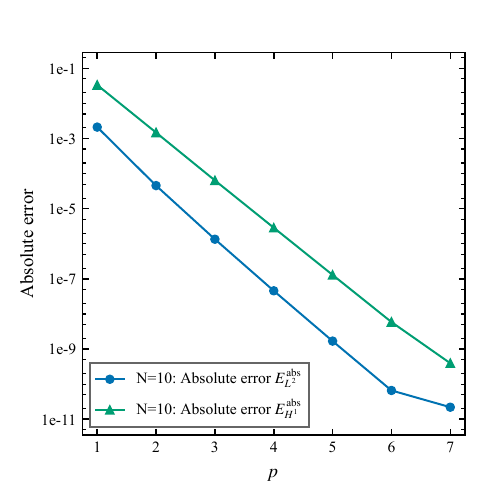}
        \caption{$N=10$.}
        \label{fig:multiple_interface_p_enrichment_N10}
    \end{subfigure}
    \hfill
    \begin{subfigure}[t]{0.45\textwidth}
        \centering
        \includegraphics[width=\linewidth]{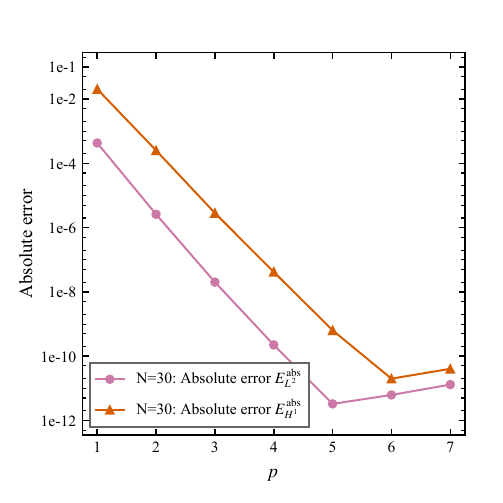}
        \caption{$N=30$.}
        \label{fig:multiple_interface_p_enrichment_N30}
    \end{subfigure}
    \caption{
        Absolute errors versus the polynomial degree $p$ for the
        multiple-interface problem.
    }
    \label{fig:multiple_interface_p_enrichment}
\end{figure*}

Figure~\ref{fig:multiple_interface_p_enrichment} shows rapid error reduction
with increasing polynomial degree until the computations approach the
$10^{-11}$--$10^{-12}$ accuracy range on the finer mesh.

\begin{figure*}[t]
    \centering
    \begin{subfigure}[t]{0.45\textwidth}
        \centering
        \includegraphics[width=\linewidth]{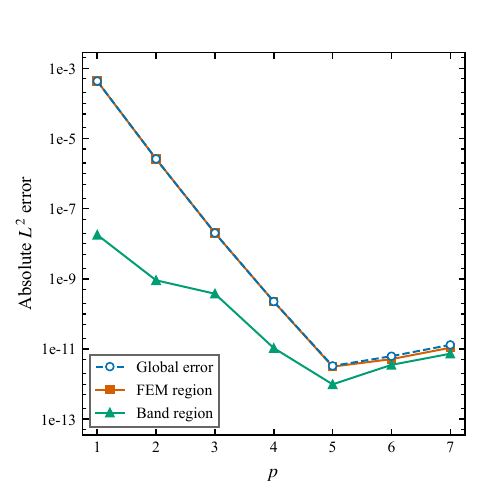}
        \caption{Absolute $L^2$ error.}
        \label{fig:multiple_interface_p_regional_l2}
    \end{subfigure}
    \hfill
    \begin{subfigure}[t]{0.45\textwidth}
        \centering
        \includegraphics[width=\linewidth]{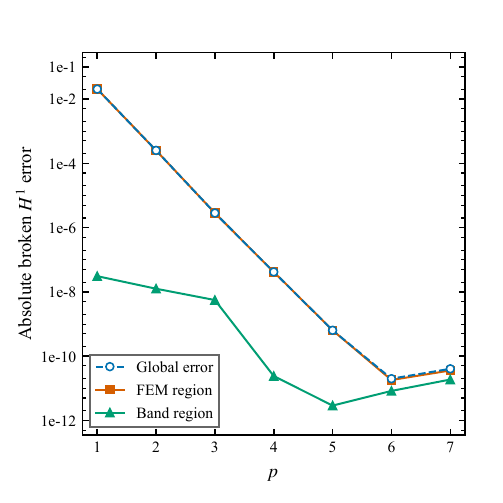}
        \caption{Absolute broken $H^1$ error.}
        \label{fig:multiple_interface_p_regional_h1}
    \end{subfigure}
    \caption{
        Global and regional errors versus $p$ for $N=30$.
    }
    \label{fig:multiple_interface_p_regional}
\end{figure*}

Before saturation, the global errors closely follow the FEM-region
errors, whereas the interface band errors remain lower. For the fixed
background mesh with $N=30$, the FEM and band regions do not change
with $p$: increasing the polynomial degree enriches the standard FEM
space and changes its coupling trace space on this fixed geometric
partition. The regional curves therefore describe how the
volume error contributions change under $p$-enrichment on a fixed
partition. Near the numerical accuracy floor, these contributions
become comparable and no longer decrease systematically. They quantify the
volume error contributions in the FEM and band regions of the coupled
solution.

\begin{figure*}[t]
    \centering
    \begin{subfigure}[t]{0.315\textwidth}
        \centering
        \includegraphics[width=\linewidth]
        {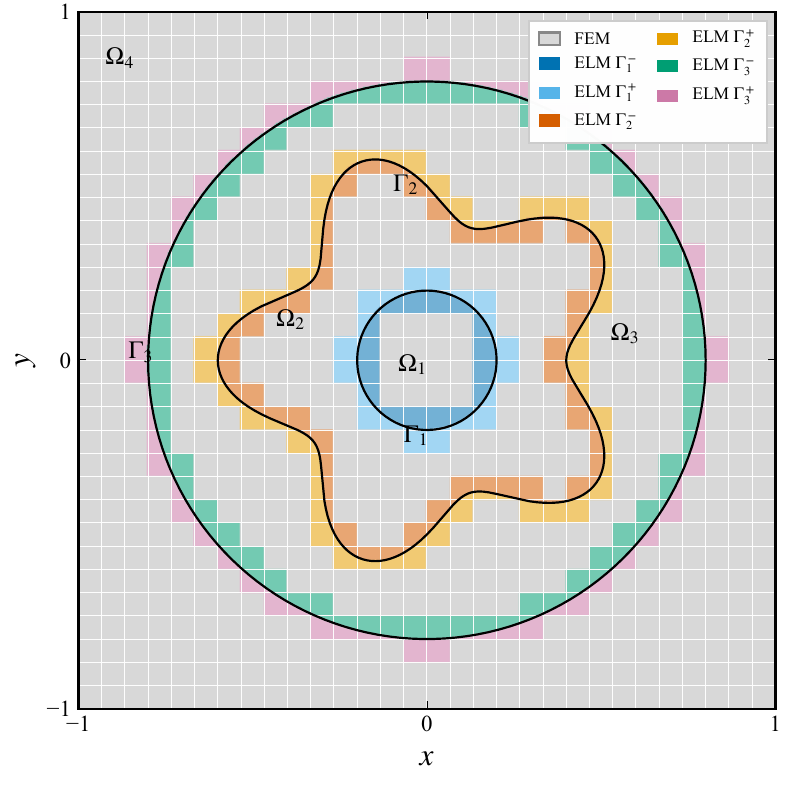}
        \caption{Hybrid FEM--RFM partition.}
        \label{fig:multiple_interface_partition}
    \end{subfigure}
    \hfill
    \begin{subfigure}[t]{0.315\textwidth}
        \centering
        \includegraphics[width=\linewidth]
        {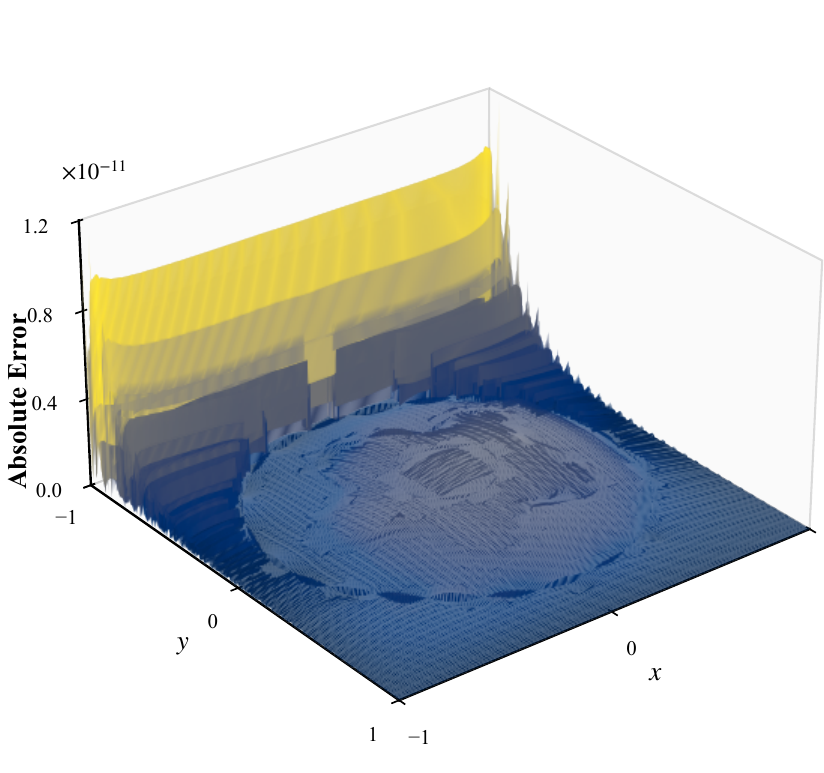}
        \caption{Pointwise solution error $|u-u_h|$.}
        \label{fig:multiple_interface_solution_distribution}
    \end{subfigure}
    \hfill
    \begin{subfigure}[t]{0.315\textwidth}
        \centering
        \includegraphics[width=\linewidth]
        {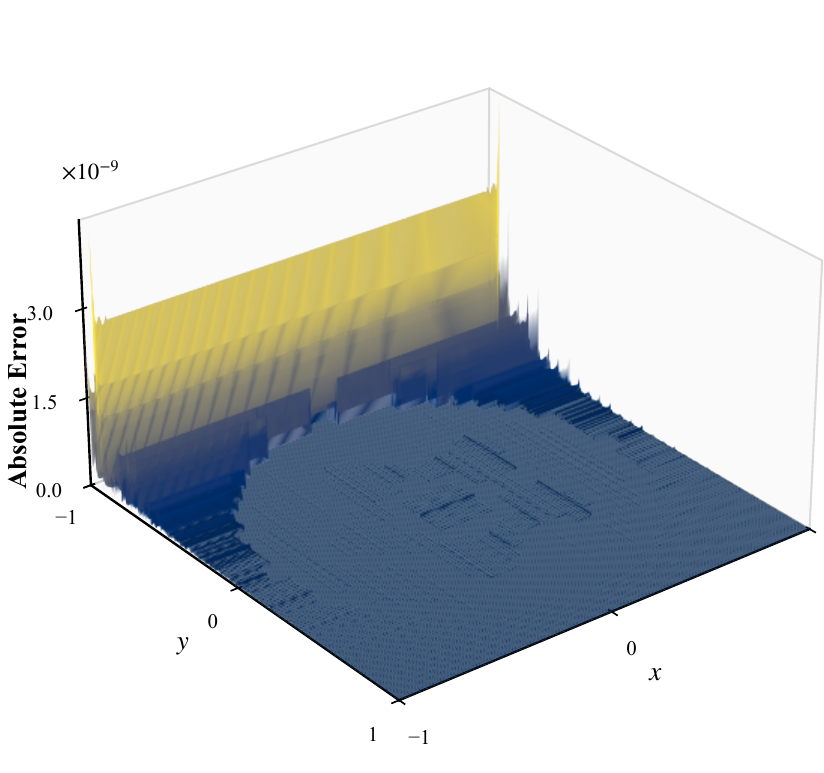}
        \caption{Pointwise gradient error
        $|\nabla u-\nabla u_h|_2$.}
        \label{fig:multiple_interface_gradient_distribution}
    \end{subfigure}
    \caption{
        Hybrid partition and pointwise absolute-error distributions for
        the $Q_5$ approximation with $N=30$. Separate two-sided local
        random-feature representations are introduced in narrow bands
        surrounding the three physical interfaces, while the remaining
        uncut cells are discretized by the finite element method.
    }
    \label{fig:multiple_interface_error_distribution}
\end{figure*}

The pointwise plots provide a representative spatial diagnostic for the
$Q_5$, $N=30$ computation. They show that the randomized approximation
remains confined to narrow neighborhoods of the three physical interfaces;
the remaining cells retain the standard finite element discretization.

\subsection{Extension diagnostics}
\label{app:extension_visual_diagnostics}

\begin{figure*}[htbp]
    \centering
    \includegraphics[width=0.90\textwidth]
    {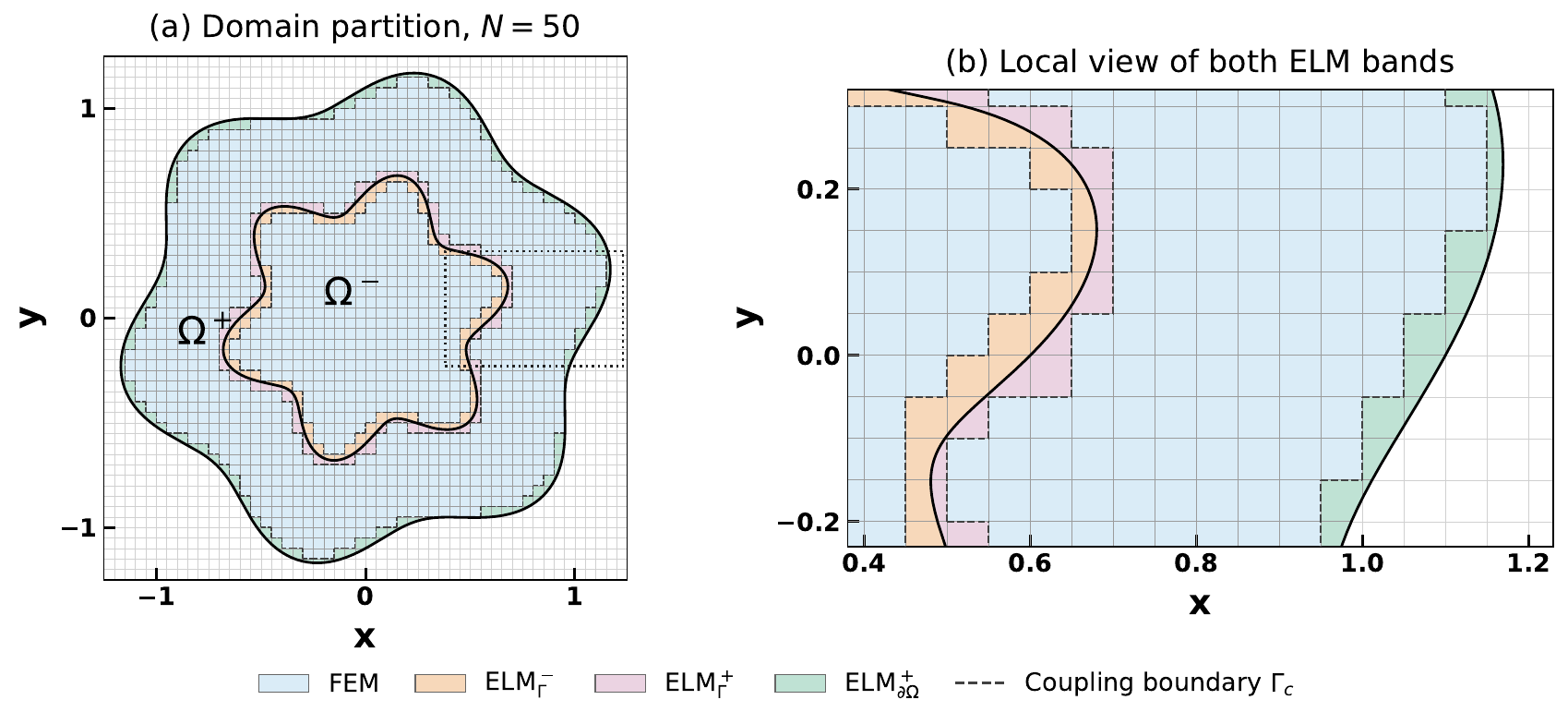}
    \caption{
        Hybrid FEM--RFM partition for the irregular-domain benchmark
        on the $N=50$ background mesh.
        Two local random-feature bands are introduced on the two sides
        of the internal interface $\Gamma$, and an additional one-sided
        band is introduced along the irregular outer boundary
        $\partial\Omega$. The remaining cells form the finite element
        region. The dashed grid-aligned curves denote the artificial
        coupling boundaries $\Gamma_c$. The right panel provides a local
        view of the internal-interface and outer-boundary bands together
        with the intervening FEM region.
    }
    \label{fig:irregular_domain_partition}
\end{figure*}

For the irregular-domain benchmark, Figure~\ref{fig:irregular_domain_partition}
shows the hybrid partition used in the computation. The additional
random-feature band is confined to cells intersected by the physical outer
boundary, while the intervening uncut cells remain in the FEM region.

\begin{figure*}[t]
    \centering
    \begin{subfigure}[t]{0.33\textwidth}
        \centering
        \includegraphics[width=\linewidth]
            {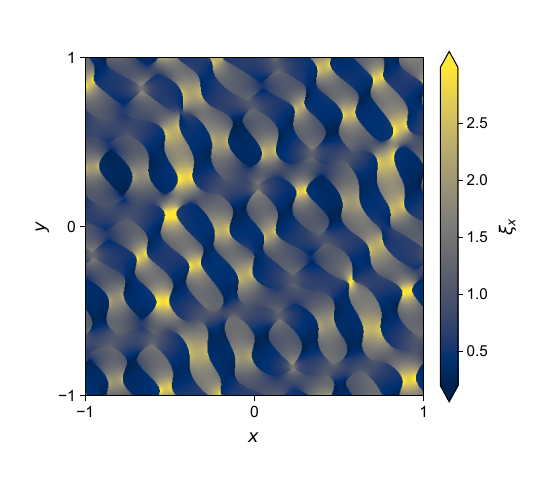}
        \caption{Reference, $N=4096$.}
    \end{subfigure}
    \hfill
    \begin{subfigure}[t]{0.33\textwidth}
        \centering
        \includegraphics[width=\linewidth]
            {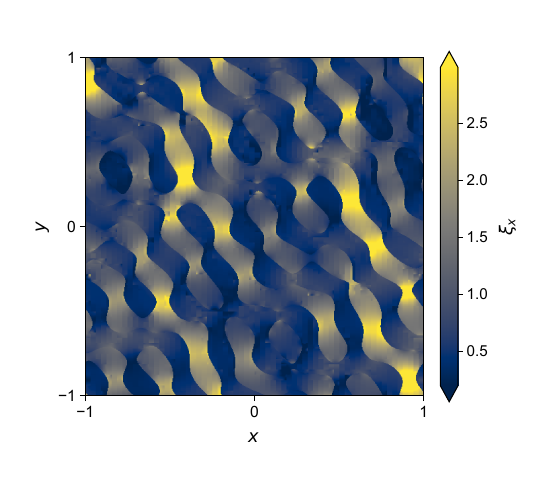}
        \caption{Hybrid $Q_1$, $N=64$.}
    \end{subfigure}
    \hfill
    \begin{subfigure}[t]{0.33\textwidth}
        \centering
        \includegraphics[width=\linewidth]
            {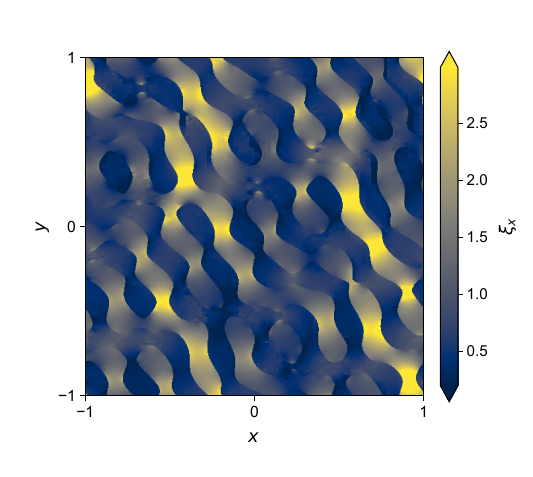}
        \caption{Hybrid $Q_2$, $N=64$.}
    \end{subfigure}
    \caption{Local negative-temperature-gradient component $\xi_x$.}
    \label{fig:complex-gradient}
\end{figure*}

For the bicontinuous composite, Figure~\ref{fig:complex-gradient} provides
the corrected thermal-gradient field complementary to the heat-flux
comparison retained in the main text. As with the heat flux, the figure is
used only as a qualitative local-field diagnostic.

\begin{figure}[htbp]
    \centering
    \begin{subfigure}[t]{0.485\textwidth}
        \centering
        \includegraphics[width=\linewidth]
        {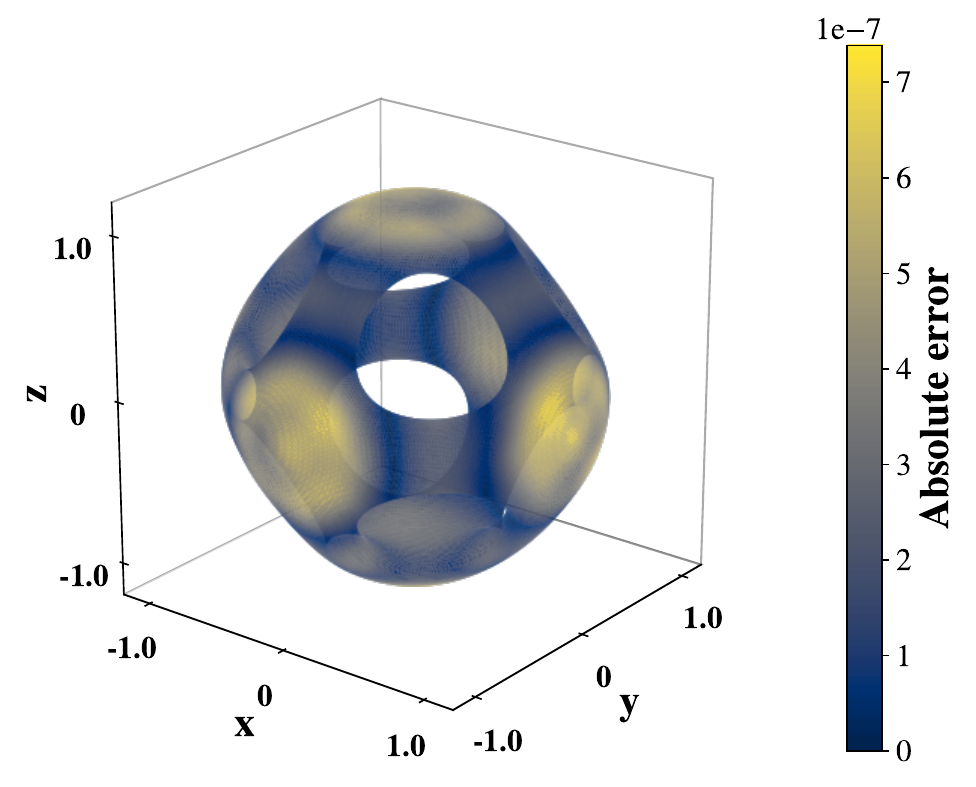}
        \caption{Case~1.}
        \label{fig:three_dimensional_case1_interface_error}
    \end{subfigure}
    \hfill
    \begin{subfigure}[t]{0.485\textwidth}
        \centering
        \includegraphics[width=\linewidth]
        {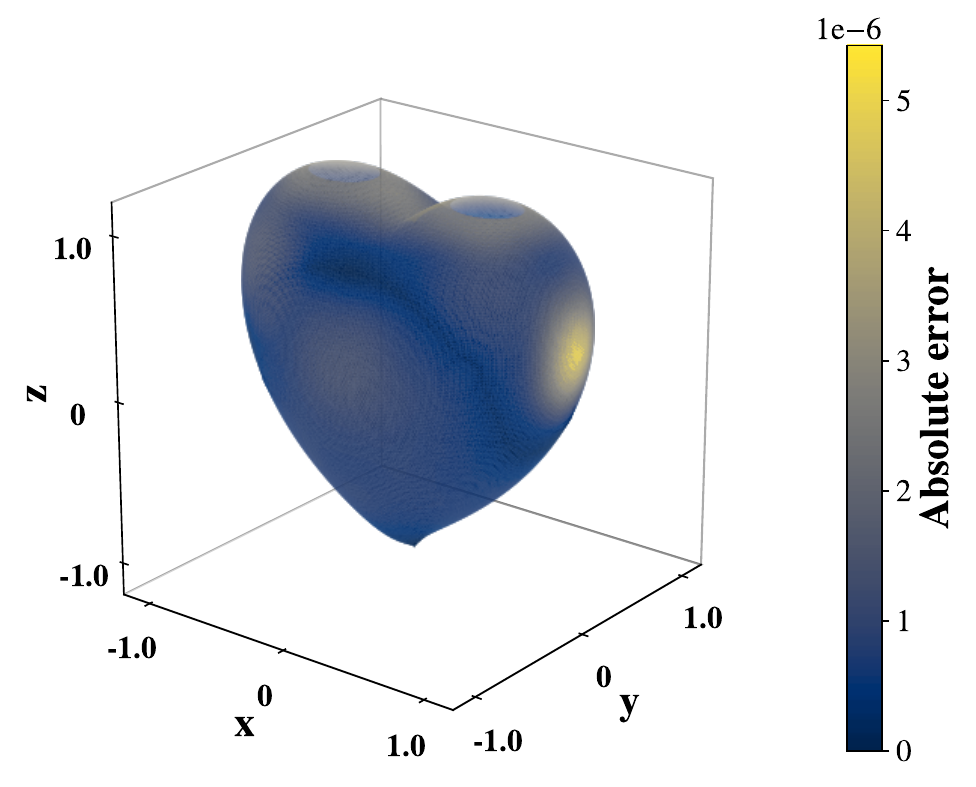}
        \caption{Case~2.}
        \label{fig:three_dimensional_case2_interface_error}
    \end{subfigure}
    \caption{
        Pointwise absolute errors of the exterior numerical trace on the
        physical interface for representative $Q_3$ computations with
        $N=20$.
    }
    \label{fig:three_dimensional_interface_errors}
\end{figure}

For the three-dimensional tests, Figure~\ref{fig:three_dimensional_interface_errors}
shows representative $Q_3$, $N=20$ pointwise errors of the exterior numerical
trace on the two interfaces. These plots supplement the $Q_1$--$Q_2$
refinement study and are not used to infer a $Q_3$ convergence rate.

\begin{figure}[htbp]
\centering
\includegraphics[width=0.9\linewidth]{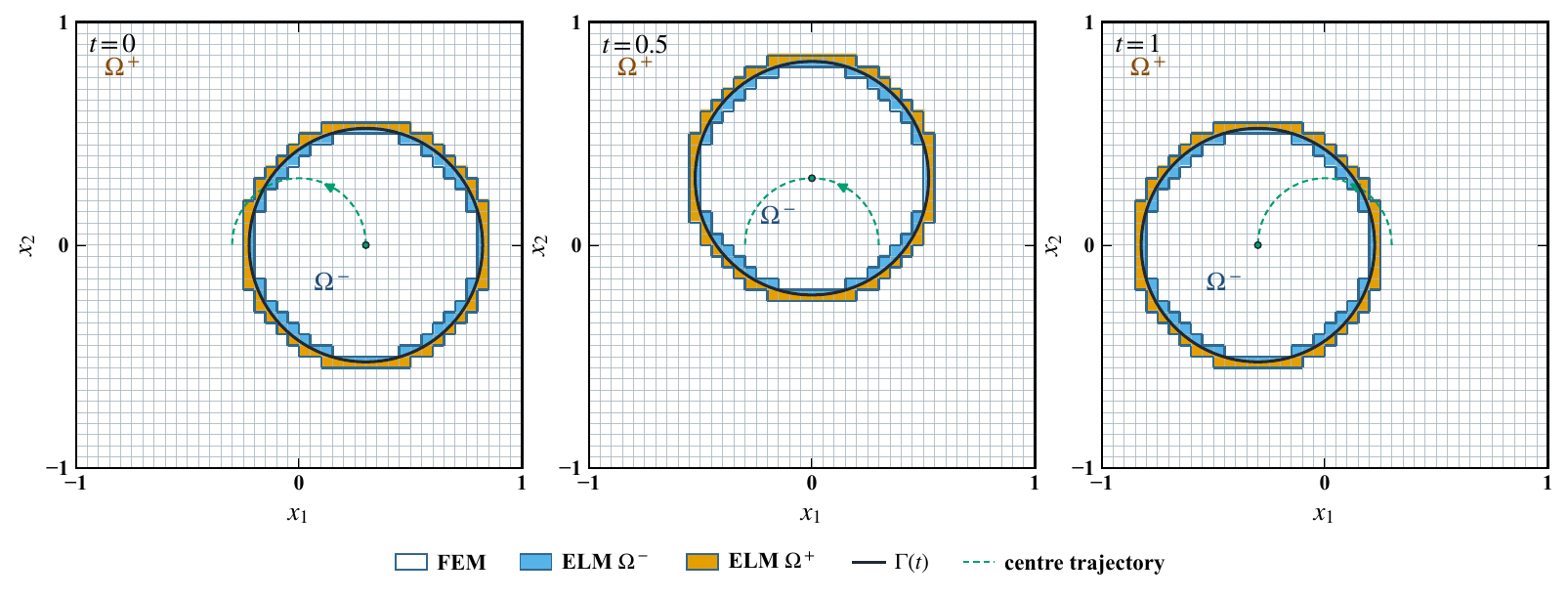}
\caption{Hybrid partition for the moving-interface problem in
Section~\ref{subsec:moving_interface_parabolic} at $N=40$.}
\label{fig:moving_interface_partition}
\end{figure}

Finally, Figure~\ref{fig:moving_interface_partition} illustrates the fixed-band
geometry used in the moving-interface experiment at $N=40$.

\end{document}